\documentclass[11pt,a4paper]{article}
\usepackage[latin1]{inputenc}
\usepackage{amsmath}
\usepackage{amsfonts}
\usepackage{amssymb}
\usepackage{graphicx}
\usepackage{algorithmic}
\usepackage{algorithm}
\usepackage{geometry}

\newtheorem{theorem}{Theorem}[section]
\newtheorem{proposition}{Proposition}[section]

\newtheorem{corollary}[theorem]{Corollary}
\newtheorem{example}[theorem]{Example}
\newtheorem{remark}[theorem]{Remark}
\def\roff  {\mbox{\boldmath$\varepsilon$}}
\def\bfeta  {\mbox{\boldmath$\eta$}}
\def\bfmu  {\mbox{\boldmath$\mu$}}
\def\bfzeta  {\mbox{\boldmath$\zeta$}}
\def\bfxi  {\mbox{\boldmath$\xi$}}
\def\bfa  {\mbox{\boldmath$a$}}
\def\bfkappa  {\mbox{\boldmath$\kappa$}}

\newcommand{\Vanderm}{\mathcal{V}}
\newcommand{\Hankel}{\mathcal{H}}

\newcommand{\spabs}{{\boldsymbol |\!\!|\!\!|}}

\numberwithin{equation}{section}

\begin{document}

\title{Numerical methods for accurate computation of the eigenvalues of Hermitian matrices and the singular values of general matrices}

\author{Zlatko Drma\v{c}${}^\dag$
	\\[2mm]
	${}^\dag$Department of Mathematics, Faculty of Science, University of Zagreb, Croatia
}
\providecommand{\keywords}[1]{\textbf{\textit{Key words:}} #1}
\date{}

\maketitle

\begin{abstract}
This paper offers a review of numerical methods for computation of the eigenvalues of Hermitian matrices and the singular values  of general and some classes of structured matrices. The focus is on the main principles behind the methods that guarantee high accuracy even in the cases that are ill-conditioned for  the conventional methods. 
First, it is shown that a particular structure of the errors in a finite precision implementation of an algorithm allows for a much better measure of sensitivity and that computation with high accuracy is possible despite a large classical condition number. Such structured errors incurred by finite precision computation are in some algorithms e.g. entry-wise or column-wise small, which is much better than the usually considered errors that are in general small only when measured in the Frobenius matrix norm. Specially tailored perturbation theory for such structured perturbations of Hermitian matrices guarantees much better bounds for the relative errors in the computed eigenvalues.  
Secondly, we review an unconventional approach to accurate computation of the singular values and eigenvalues of some  notoriously ill-conditioned structured matrices, such as e.g. Cauchy, Vandermonde and Hankel matrices. The distinctive feature of accurate algorithms is using the intrinsic parameters that define such matrices to obtain a non-orthogonal factorization, such as the \textsf{LDU} factorization, and then computing the singular values of the product of thus computed factors. The state of the art software is discussed as well.
\end{abstract}	

\keywords backward error, condition number, eigenvalues, Hermitian matrices, Jacobi method, LAPACK, perturbation theory, rank revealing decomposition, singular value decomposition

\section{Introduction}

  In real world applications, numerical computation is done with errors (model errors, measurement errors, linearization errors, truncation/discretization errors, finite computer arithmetic errors). This calls for caution when interpreting the computed results. For instance, any property or function value we obtain from finite precision computation with a nontrivial matrix $A$ stored in the computer memory (for instance, the rank or the eigenvalues of $A$) very likely holds true for some unknown $A+\delta A$ in the vicinity of $A$, but not for $A$. In order to estimate the level of accuracy that can be expected in the output, we need to know the level of initial uncertainty in the data, the analytical properties of the function of $A$ that we are attempting to compute,
  the numerical properties of the algorithm used and the parameters of the computer arithmetic.

  A better understanding of the sensitivity of numerical problems, together with the adoption of
  new paradigms in the algorithmic development over the last few decades have opened new possibilities, allowing for high accuracy solutions to problems that were previously considered numerically intractable. In this paper we give an overview of such advances as regards the computation to high accuracy of the eigenvalues  of Hermitian matrices and the singular values  of general and some special classes of matrices. The focus is on the main principles, and technical details will be mostly avoided.

  Computing the eigenvalues with high accuracy means that
  for each eigenvalue (including the tiniest ones, much smaller than the norm of the matrix) as many correct digits are computed as warranted by the data. In other words, for the eigenvalues
  $\lambda_1\geq\cdots\geq\lambda_n$ of a nonsingular Hermitian matrix $H=H^*\in\mathbb{C}^{n\times n}$ and their computed approximations
  $\widetilde\lambda_1\geq\cdots\geq\widetilde\lambda_n$
  we want a bound of the form
  \begin{equation}\label{zd:eq:intro:rel_error}
  \max_{i=1:n} \frac{|\widetilde\lambda_i-\lambda_i|}{|\lambda_i|} \leq \bfkappa \cdot O(\roff) ,
  \end{equation}
  where $\bfkappa$ represents a hopefully  moderate condition number, and $\roff$ is the round-off unit of the computer arithmetic.

  For this kind of accuracy, the standard paradigm (algorithm based on
  orthogonal transformations, small backward error and perfect stability of the symmetric eigenvalue problem)
  is not good enough. Namely, the conventional approach of showing that the computed $\widetilde\lambda_i$'s are the exact eigenvalues of a nearby
  $H+\delta H$ with $\|\delta H\|_2 \leq O(\roff) \|H\|_2$, and then applying Weyl's theorem, which guarantees
  that $\max_{i=1:n} |\widetilde\lambda_i-\lambda_i| \leq \|\delta H\|_2$, yields, for each eigenvalue index $i$,
  \begin{equation}\label{zd:eq:intro:rel_norm_error}
  \frac{|\widetilde\lambda_i-\lambda_i|}{\|H\|_2} \leq O(\roff),\;\;\mbox{i.e.}\;\;
  \frac{|\widetilde\lambda_i-\lambda_i|}{|\lambda_i|} \leq O(\roff) \frac{\|H\|_2}{|\lambda_i|} \leq O(\roff) \frac{\|H\|_2}{|\lambda_n|} = O(\roff) \|H\|_2 \|H^{-1}\|_2.
  \end{equation}
  Clearly, (\ref{zd:eq:intro:rel_norm_error}) will give a satisfactory bound of the form (\ref{zd:eq:intro:rel_error})
  only for absolutely large eigenvalues (those with $|\lambda_i|$ of the order of the norm $\|H\|_2$), while the relative error in the smallest eigenvalues ($|\lambda_i|\ll \|H\|_2$)
  is up to $O(\roff) \kappa_2(H)$, where $H$ is assumed nonsingular, $\kappa_2(H)=\|H\|_2 \|H^{-1}\|_2$ is the condition number, and $\|\cdot\|_2$ is the spectral operator norm, induced by the Euclidean vector norm.

  Hence, in the conventional setting, the eigenvalues of Hermitian/symmetric matrices are not always perfectly well conditioned in the sense that we can compute them in finite precision with small relative error (\ref{zd:eq:intro:rel_error}). We need to identify classes of matrices that allow for such high relative accuracy. To that end, we may need to restrict the classes of permissible perturbations -- instead of in matrix norm, one may consider finer, entry-wise changes in matrix entries. As a result of such stronger requirements of relative accuracy, there will be a new condition number able to distinguish between well- and ill-behaved matrices with  respect to such perturbations. Hence, for some classes of matrices we will
  	be able to compute even the tiniest eigenvalues  even if $\kappa_2(H)$ is extremely large. For that, however, we will have to rethink and redefine the paradigms of algorithm development. 
  	
  	For the sake of brevity, in this review we do not discuss the accuracy of the computed eigenvectors and the singular vectors. This is an important issue, and interested readers will find the relevant results in the provided references.
  	
  	The new structure of the perturbation (finer than usually required by $\|\delta H\|_2 \leq O(\roff) \|H\|_2$)  and the
  	condition number governing high relative accuracy are not invariant under general orthogonal
  	similarities. This means that using an algorithm based on orthogonal transformations does not automatically
  	guarantee results that are accurate in the sense of (\ref{zd:eq:intro:rel_error}). Some algorithms are more accurate than the others, see \cite{dem-ves-92}. For best results, separate perturbation theories and numerical algorithms have to be developed for the positive definite and the indefinite matrices. 
  	
  	Analogous comments apply to the computation of the singular values $\sigma_1\geq\cdots\geq\sigma_n$ of $A\in\mathbb{C}^{m\times n}$ -- conventional algorithms in general cannot approximate a small singular value $\sigma_i$ to any correct digit  if $\sigma_i < \roff \sigma_1$, despite the fact that only orthogonal or unitary transformations are used. 

This review of the development of new theory and new algorithms is organized as follows. For the readers' convenience, in \S \ref{S=back-stability} we first give a brief review of the key notions of backward stability, perturbation theory, condition number and forward error. 
Then,  in \S \ref{S=HPD}, we review the state of the art numerical methods for  computing the eigenvalues of real symmetric and Hermitian matrices. Brief description of the algorithms in \S \ref{ss=Classical-Methods} is followed by a general framework for the classical backward error analysis in \S \ref{SS=backward-stability}, and its 
limitations with respect to the high accuracy of the form (\ref{zd:eq:intro:rel_error}) is shown in \S \ref{zd:SSS:Example_3x3} using a $3\times 3$ symmetric matrix as a case study. The conditions for achieving (\ref{zd:eq:intro:rel_error}) are analyzed in \S \ref{SS=posdef-accurate} for positive definite matrices, and in \S \ref{SSS=symm-jac-pd} we show that the symmetric Jacobi algorithm in finite precision arithmetic satisfies these conditions, which makes it provably  more accurate than any tridiagonalization based algorithm \cite{dem-ves-92}.
This theory does not include indefinite matrices, which are analyzed separately in \S \ref{S=HID}.  It will become clear that there is a fundamental difference between the two classes.  
We conclude the numerical computation with positive definite matrices by arguing in \S \ref{SS=implicit-pd} that in many cases such matrices are best given implicitly by a factor $A$ such that $A^*A = H$, and that accurate eigenvalue computation follows from accurate computation of the singular values of $A$.

In \S \ref{S=SVD} we study algorithms for computing the singular values to high relative accuracy.
After reviewing the bidiagonalization based methods in \S \ref{SS=bidiagSVD}, and the one-sided Jacobi \textsf{SVD} in \S  \ref{SS=One-sided-jacobi} and its preconditioned version in \S \ref{SS=J+QRCP}, 
in \S \ref{SSS=Impl-Jacobi-eig},  we show how the combination of the Cholesky factorization and the one-sided Jacobi \textsf{SVD} computes the eigenvalues of general (non-structured)  positive definite matrices to the optimal accuracy (\ref{zd:eq:intro:rel_error}) permitted by the perturbation theory.
 Section \ref{S=PSVD} reviews accurate computation of the \textsf{SVD} of certain products of matrices (\textsf{PSVD}), which is the core procedure for the new generation of highly accurate \textsf{SVD}  algorithms, based on the so-called rank-revealing decompositions (\textsf{RRD}s). In \S \ref{S=RRD+PSVD} we illustrate the \textsf{RRD}+\textsf{PSVD} concept in action.
In particular, we discuss structured matrices such as the Cauchy, Vandermonde and Hankel matrices, which  are the key
objects in many areas of numerical mathematics, in particular in rational approximation
theory \cite{Gutknecht:1982:RPC}, \cite{Gonnet:2011:RRI}, \cite{haut-beylkin-coneig-2011}, where e.g. the coefficients of the
approximant are taken from the singular vectors corresponding to small singular values of
certain matrices of these kinds. The fact that these structured matrices can be extremely ill-conditioned
is the main obstacle that precludes turning powerful theoretical results into practical numerical procedures. We review recent results that allow for highly accurate computations even in extremely ill-conditioned cases. 

Section \ref{S=HID} is devoted to accurate computation of the eigenvalues of Hermitian indefinite matrices.
The key steps towards understanding the sensitivity of the eigenvalues  are reviewed in \S \ref{SS=HID-pert-theory}, and in \S \ref{SS=sym-indef-fact}, \S \ref{SS=OJ} we review numerical algorithms that compute the eigenvalues of indefinite matrices to the accuracy deemed possible by the corresponding perturbation theory. Three different approaches are presented, and, interestingly, all based on the Jacobi algorithm but with some nonstandard features. In \S \ref{SS=JJ}, the problem is transformed to a generalized eigenvalue problem,  and the Jacobi diagonalization process is executed using transformations that are orthogonal in an indefinite inner product   whose signature is given by the inertia of $H$.  In \S \ref{SSS=ISJM}, the classical Jacobi algorithm is carefully implemented implicitly on an accurately computed symmetric indefinite factorization, and in \S \ref{SS=OJ} the spectral decomposition of $H$ is carefully extracted from its accurate SVD.  Extending the results to classes of non-symmetric matrices is a challenging problem and in \S \ref{SS=nonsymm-TN} we briefly review the first results in that direction.

\section{Backward stability, perturbation theory and condition number}\label{S=back-stability}
The fundamental idea of backward stability is introduced by Wilkinson \cite{Wilinson-RoundingErrors-63}, \cite{Wilkinson-AEP-65}.
In an abstract formulation, we
want to compute $Y=\mathcal{F}(X)$ using an algorithm $\mathcal{A}_{\mathcal F}(X)$ that
returns only an approximation $\widetilde{Y}$ of $Y$. In the \emph{backward error
	analysis} of the computational process $\mathcal{A}_{\mathcal F}(X)$, we prove existence of a small perturbation $\delta X$ of the input data $X$ such that
$\widetilde{Y}=\mathcal{F}(X+\delta X)$. If the perturbation $\delta X$, called \emph{backward error}, is acceptably small relative to $X$ (for instance, of the same order as
the initial uncertainty $\delta_0 X$  already present in $X=X_{\textrm{true}}+\delta_0 X$, where $X_{\textrm{true}}$ is the unaccessible exact value), the computation
of $\widetilde{Y}$ by $\mathcal{A}_{\mathcal F}(\cdot)$ is considered \emph{backward stable}.

 The sources of error can be the use of finite precision arithmetic or any other approximation scheme. Sometimes, a backward error is constructed artificially just in order to justify the computed output. 
 For instance, if we compute an approximate eigenvalue $\lambda$, with the corresponding eigenvector $v\neq \mathbf{0}$, of the matrix $X$, the residual $r = Xv-\lambda v$ will be in general nonzero, but hopefully small. We can easily check that $( X + \delta X) v = \lambda v$, with $\delta X = - r v^*/(v^* v)$, i.e. we have found an exact eigenpair $\lambda, v$ of a nearby matrix $X+\delta X$, with $\|\delta X\|_2=\|r\|_2/\|v\|_2$. If, for given $v$, we choose $\lambda$ as the Rayleigh quotient $\lambda=v^* X v/v^*v$, then $r^*v=0$ and $(X +\Delta X)v=\lambda v$ with Hermitian $\Delta X=\delta X + (\delta X)^*$, which is favorable interpretation if $X$ is Hermitian: we have solved a nearby Hermitian problem.

Backward stability does not automatically imply that $\widetilde{Y}$
is close to $Y$. It is possible that a backward stable algorithm gives utterly wrong results even in the case when the computation is considered stable; see \S \ref{zd:SSS:Example_3x3} below for an example and its discussion. The error in the result (\emph{forward error}) $\delta Y = \widetilde{Y}-Y = \mathcal{F}(X+\delta X) - \mathcal{F}(X)$
will depend on the function $\mathcal{F}(\cdot)$,  i.e. on its sensitivity to the change in the argument $X$, and on the size and the structure of $\delta X$.
This sensitivity issue is the subject of \emph{perturbation theory}.
If the forward error $\delta Y$ is small, the algorithm is called \emph{forward stable}.
Claiming backward stability of an algorithm depends on how the size of the backward error is
measured, as well as on other factors, such as the structure of the backward error. For instance, if $X$ is a symmetric matrix, it is desirable to prove existence of a symmetric perturbation $\delta X$, see e.g. \cite{Smoktunowicz:1995:NSC}.  If $X$ lives in a normed space $(\mathcal{X},\|\cdot\|_x)$, then
we usually seek a bound of the form $\|\delta X\|_x \leq \epsilon \|X\|_x$. The error in the computed result, which is assumed to live in a normed space $(\mathcal{Y},\|\cdot\|_y)$,  is estimated as $\|\delta Y\|_y \leq C \epsilon \|Y\|_y$. The amplification factor $C$ is  the \emph{condition number}.

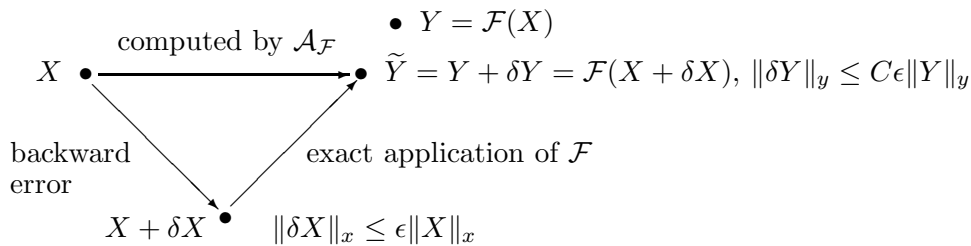
\begin{figure}[ht]
	\begin{center}
		\setlength{\unitlength}{5ex} \centering
		\begin{picture}(10,4.0)(1.5,0)
		\put(-0.10,2.5){$X\;\;\bullet$} \put(4.95,2.5){{$\bullet\;\;\widetilde{Y}=Y+\delta Y=\mathcal{F}(X+\delta X)$}, $\|\delta Y\|_y \leq C \epsilon \|Y\|_y$}
		\put(0.8,2.45){\vector(1,-1){2.0}} \put(3,0.5){\vector(1,1){2.0}}
		\put(0.90,2.60){\vector(1,0){3.98}} \put(-0.5,1.2){backward}
		\put(-0.5,0.7){error} 
		\put(4.2,1.2){exact application of $\mathcal{F}$}
		\put(1.2,3.0){computed by $\mathcal{A}_{\mathcal{F}}$} 
		\put(5.5,3.3){{$\bullet\;\; Y=\mathcal{F}(X)$}}
		\put(1.0,0){$X+\delta X\;\;$ \hspace{0.5cm}$\|\delta X\|_x \leq \epsilon
			\|X\|_x$}
		\put(2.8,0.2){$\bullet$}
		\end{picture}
		\caption{\label{backdiag} Commutative diagram for backward stable computation. The value $\widetilde{Y}$ computed by the algorithm $\mathcal{A}_{\mathcal{F}}(X)$ is the exact value of the function $\mathcal{F}$ at $X+\delta X$. The condition number $C$ governs the size of the forward error $\delta Y$.}
	\end{center}
\end{figure}

An abstract theory of the condition number, ill-conditioning and related problems in the setting of normed
manifolds is given e.g. in \cite{Rice-TheoryofCondition-66}. As an illustration of such an abstract analytical treatment,
we cite one result:
\begin{theorem} (Rice \cite{Rice-TheoryofCondition-66})
	Let $(\mathcal{X},\|\cdot\|_x)$, $(\mathcal{Y},\|\cdot\|_y)$ be normed linear spaces
	and $\mathcal{F}: \mathcal{X} \longrightarrow \mathcal{Y}$
	be a differentiable function. The absolute and the relative condition numbers,
	respectively, of $\mathcal{F}$ at $X_0$ are defined as
	\begin{eqnarray*}
		\alpha(\mathcal{F},X_0;\delta) &=& \inf\{C\geq 0\; :\; \|X-X_0\|_x <
		\delta \Longrightarrow \| \mathcal{F}(X)-\mathcal{F}(X_0)\|_y < C \delta \}\\
		\rho(\mathcal{F},X_0;\delta) &=& \inf\{C\geq 0\; :\; \|X-X_0\|_x < \delta \|X_0\|_x \Longrightarrow
		\| \mathcal{F}(X)-\mathcal{F}(X_0)\|_y <
		C \delta \|\mathcal{F}(X_0)\|_y\}.
	\end{eqnarray*}
	Let the corresponding asymptotic condition numbers be defined as
	\begin{displaymath}
	\alpha(\mathcal{F},X_0) = \lim_{\delta\rightarrow 0} \alpha(\mathcal{F},X_0;\delta),\;\;\;\;
	\rho(\mathcal{F},X_0) = \lim_{\delta\rightarrow 0} \rho(\mathcal{F},X_0;\delta).
	\end{displaymath}
	If $\mathcal{J}$ is the Jacobian of $\mathcal{F}$ at $X_0$, then
	$
	\alpha(\mathcal{F},X_0) = \| \mathcal{J}\|_{xy}$,
	$\rho(\mathcal{F},X_0) = \frac{\| \mathcal{J}\|_{xy}}{\|\mathcal{F}(X_0)\|_y}\|X_0\|_x,
	$
	where $\|\cdot\|_{xy}$ denotes the induced operator norm.
\end{theorem}
For a systematic study of the condition numbers, we refer to \cite{book-condition}, and for various techniques of backward error analysis see \cite{book-higham}. The backward error analysis is often represented using commutative diagrams such as in Figure \ref{backdiag}.

\subsection{Scaling and numerical stability}\label{SS:scaling-Hankel-SVD}
When doing numerical calculations and estimating errors, computing residuals, or making decisions about the ranks of matrices by
declaring certain quantities sufficiently small to be considered negligible, it is often forgotten  that those numbers represent
physical quantities in a particularly chosen system of units. A small entry in the matrix may be merely noise, but it
could also be a relevant small physical parameter. Very often, the system under consideration represents couplings
between quantities of different physical natures, each one given in its own units on a particular scale. In fact, it is
a matter of engineering design and ingenuity  to choose the units so that the mathematical model faithfully represents the
physical reality, and that the results of the computations can be meaningfully measured in appropriate norms and interpreted and used with
confidence in applications.

As an illustration, we briefly discuss one simple example: consider a state space realization of a linear time invariant (\textsf{LTI}) dynamical system
\begin{eqnarray}
\dot x(t) &=&  A x(t) + B u(t),\; x(0) = x_0 , \label{eq:LTI1}\\
y(t) &=& C x(t) .  \label{eq:LTI2}
\end{eqnarray}
In general, switching to different units implies a change of variables  \,$x(t) = S\hat{x}(t)$, where $S$
denotes the corresponding diagonal scaling matrix.
In the new set of state variables defined by \,$x(t)= S \hat x(t)$, the system goes over into
$
\dot{\hat{x}}(t) = (S^{-1}AS) \hat{x}(t) + (S^{-1}B) u(t)$, 
$y(t) = (CS) \hat{x}(t) . 
$
Hence, the state space description $(A,B,C)$ changes into
the equivalent triplet $(S^{-1} A S, S^{-1}B, CS)$.
From the system-theoretical point of view, nothing has changed, the new state space representation represents
the same physical reality: the transfer function $\mathcal{G}(s) = C(sI-A)^{-1}B = (CS) (sI - S^{-1}A S)^{-1}(S^{-1}B)$ is the same, the system poles (the eigenvalues of $A$, i.e. of $S^{-1}AS$) are also the same.

Unfortunately, this invariance is not inherited in finite precision computation. For example, the important property of stability requires the eigenvalues of the system matrix $A$ to be in the open left complex half
plane. If the numerically computed eigenvalues do satisfy that condition but are too close to the imaginary axis, how can we be certain that the system is stable?
Since the system matrix is determined up to a similarity, how can we be sure that our numerical algorithm will not be influenced
by a particular representation?

Important properties of the system (\ref{eq:LTI1},\ref{eq:LTI2}), such as controllability and observability, are encoded in the symmetric positive semidefinite matrices \,$H= \int_0^\infty e^{t A} BB^T e^{t A^T}dt$\, and \,$M= \int_0^\infty
e^{t A^T}C^T C e^{t A} dt$, called system Gramians,
which are computed as the solution of the dual pair of Lyapunov equations
$
{{A}H + H {A}^{T} = - {B} {B}^{T}}, \;\;
{{A}^{T}M+M{A}  = - {C}^{T}{C}}.
$
The joint spectral properties of the Gramians provide information on deep structural properties of
the system.
The key quantities in this respect are the Hankel singular values, defined as the square roots \,$\sigma_i = \sqrt{\lambda_i(HM)}$\, of the eigenvalues of $HM$. These are proper invariants of the system, independent of the
state space realization (\ref{eq:LTI1},\ref{eq:LTI2}).
It can be easily checked that changing into $\hat{x}(t)$ changes system Gramians by the so called contragredient transformation
$H \longrightarrow \widehat{H} = S^{-1} H S^{-T}$, $M \longrightarrow \widehat{M} = S^T M S$, which does not change the Hankel singular values, since $\widehat{H}\widehat{M} = S^{-1}(HM)S$.

The numerics, on the other hand,  may react sharply. A change of units (scaling) changes classical condition numbers
$\kappa_2(A)$, $\kappa_2(H)$, $\kappa_2(M)$ thus potentially making an algorithm numerically inaccurate/unstable, while, at the same time, the underlying problem is the same. With a particularly
chosen $S$ we can manipulate the numerical rank of any of the two Gramians, and thus mislead numerical algorithms.
Is this acceptable?
If a rank decision has to be made, and if the determined \emph{numerical rank} (cf. \cite{Golub-Klema-Stewart-Numerical_rank}) sharply changes with the change of
physical units in which the variables are expressed, one definitely has to ask many nontrivial questions.
It is also possible that two algebraically equivalent methods for checking controllability of a LTI system give
completely different estimates of numerical rank.
For an excellent discussion on these topics we refer to \cite{Paige_Numerics_in-Control}.

\section{Computing eigenvalues of  Hermitian matrices}\label{S=HPD}
Numerical computation of the eigenvalues and eigenvectors of Hermitian matrices is considered as an example of a perfect computational process. This is due to several important spectral properties of Hermitian matrices, see e.g. \cite[Ch. 4]{hor-joh-90}.

The Schur form of Hermitian matrices is diagonal: for any Hermitian $H\in\mathbb{C}^{n\times n}$ there exists a unitary $U\in\mathbb{C}^{n\times n}$ and a real diagonal $\Lambda=\mathrm{diag}(\lambda_i)_{i=1}^n$ such that $H = U \Lambda U^*$. If $U=\begin{pmatrix} u_1 & u_2 & ,\ldots,& u_n\end{pmatrix}$ is the column partition of $U$, then $H u_i = \lambda_i u_i$, $i=1,\ldots, n$. Hence, the diagonalization is performed by a unitary (or real orthogonal if $H$ is real symmetric) matrix of eigenvectors, which allows for numerical algorithms based only on unitary transformations (unitary transformations are preferred in finite precision arithmetic because they preserve relevant matrix norms, e.g., $\|\cdot\|_2$ and the Frobenius norm $\|\cdot\|_F$, and thus will not inflate initial data uncertainties and unavoidable rounding errors).

Furthermore, the eigenvalues of $H\equiv H^*$ have far-reaching variational characterizations.
\begin{theorem}\label{TM-minmax}
Let $H$ be $n\times n$ Hermitian with eigenvalues 	
$\lambda_1\geq\lambda_2\geq\cdots\geq\lambda_n$.
Then
$$
\lambda_j = \max_{{\mathcal{S}_j}} \min_{x\in\mathcal{S}_j\setminus \{0\}} \frac{x^* H x}{x^* x}
$$	
where the maximum is taken over all $j$-dimensional subspaces $\mathcal{S}_j$ of $\mathbb{C}^n$.
\end{theorem}
This characterization generates a sound perturbation theory that provides a basis for assessing the accuracy of numerical methods. For the sake of completeness  and for the reader's convenience, we cite one of the Weyl-type theorems:
\begin{theorem}\label{zd:TM:Weyl} (Weyl)
	Let $H$ and $H+\delta H$ be Hermitian  matrices with eigenvalues
	$\lambda_1\geq\lambda_2\geq\cdots\geq\lambda_n$ and $\widetilde{\lambda}_1\geq\widetilde{\lambda}_2\geq \cdots\geq\widetilde{\lambda}_n$,
	respectively. Write $\widetilde{\lambda}_i=\lambda_i+\delta\lambda_i$. Then
	$
	\max_{i=1:n} |\delta\lambda_i| \leq \|\delta H\|_2.
	$
\end{theorem}
For a {more} detailed overview {of the specific spectral properties of Hermitian matrices, including the perturbation theory,} we refer to \cite{ste-sun-90}, \cite{Bhatia-MatrixAnalysis-1997}.
\subsection{Classical methods}\label{ss=Classical-Methods}
The common paradigm of modern numerical algorithms for computing a unitary eigenvector matrix $U$ and the real diagonal $\Lambda$ is to build a sequence of unitary similarities such that
\begin{equation}\label{eq:Hk-Lambda}
H^{(k+1)} = (U^{(k)})^* \cdots ((U^{(2)})^* ((U^{(1)})^* H U^{(1)}) U^{(2)}) \cdots U^{(k)} \longrightarrow \Lambda=
\left(\begin{smallmatrix} \lambda_1 & 0& 0\cr
0 & \ddots & 0 \cr                                                             0  & 0     & \lambda_n\end{smallmatrix}\right),\;\;\quad \mbox{as }\ \  k\longrightarrow\infty.
\end{equation}
The accumulated infinite product $U^{(1)} U^{(2)}\cdots U^{(k)}\cdots$ provides information about the eigenvectors and eigenspaces (in case of
multiple eigenvalues).

The choice of unitary matrices $U^{(i)}$ defines the specific algorithm. For a detailed overview with references we recommend \cite{parlett-sevp-98},  \cite[Ch. 8]{golub-vnl-4} and \cite[Ch. 55]{hogben14}.
The two typical classes of methods are:
\begin{itemize}
	\item \emph{Tridiagonalization-based methods}: The matrix $H$ is first reduced to a Hermitian tridiagonal matrix $T$:
	\begin{equation}\label{eq:tridiagonal}
	V^* H V = T = \left(\begin{smallmatrix} \alpha_1 & \beta_1  &        &   \cr
	\beta_1  & \alpha_2 & \ddots &   \cr
	& \ddots   & \ddots & \beta_{n-1} \cr
	&          & \beta_{n-1} & \alpha_n \end{smallmatrix}\right),
	\end{equation}
	where $V$ denotes a unitary matrix composed as the product of $n-2$ Householder reflectors.
	The tridiagonalization can be illustrated in the $4\times 4$ case as
	$$
	H^{(1)}=V_1^* H V_1 = \left(\begin{smallmatrix}
	\star & \star & 0 & 0\cr
	\star & \star & \times & \times\cr
	0 & \times & \times & \times\cr
	0 & \times & \times & \times\end{smallmatrix}\right),\;\;
	H^{(2)}=V_2^* H^{(1)} V_2 = \left(\begin{smallmatrix}
	\star & \star & 0 & 0\cr
	\star & \star & \star & 0\cr
	0 & \star & \star & \star\cr
	0 & 0 & \star & \star\end{smallmatrix}\right) = T = (V_2^* V_1^*) H (V_1 V_2),
	$$		
    {where each $\star$  denotes an entry which has already
    been modified by the algorithm and set to its final value}.

	In the second stage,
	fast algorithms specially tailored for {tridiagonal} matrices, such as \textsf{QR}, divide and conquer, \textsf{DQDS},
    inverse iteration or the \textsf{MRRR} method, are deployed to compute the spectral decomposition of $T$ as
	$T = W \Lambda W^*$. Assembling back, the spectral decomposition of $H$ is obtained as  $H = U\Lambda U^*$ with $U = V W$. For studying excellent tridiagonal eigensolvers with many fine mathematical and numerical details we recommend \cite{mgu-eis-tridiag-95}, \cite{PARLETT2000121}, \cite{dhillon-parlett-2004}, \cite{DHILLON20041}, \cite{MRRR-2006}. 
	\item \emph{Jacobi type methods}:
	The classical Jacobi method
	generates a sequence of unitary congruences, ${H}^{(k+1)} =
	({U}^{(k)})^* {H}^{(k)} {U}^{(k)}$, where
	${U}^{(k)}$
	differs from the identity only at some cleverly chosen positions
	$(i_k,i_k)$, $(i_k,j_k)$, $(j_k,i_k)$, $(j_k,j_k)$, with
	$$
	\begin{pmatrix}{U}^{(k)}_{i_k,i_k} &
	{U}^{(k)}_{i_k,j_k}\cr\\[0.1pt] {U}^{(k)}_{j_k,i_k} &
	{U}^{(k)}_{j_k,j_k}\end{pmatrix} =
	\begin{pmatrix}\cos\phi_k & e^{i\psi_k}\sin\phi_k \cr
	-e^{-i\psi_k}\sin\phi_k & \cos\phi_k\end{pmatrix}.
	$$
	The
	angles $\phi_k$, $\psi_k$ of the $k$-th transformation are determined to annihilate the
	$(i_k,j_k)$ and $(j_k,i_k)$ positions in ${H}^{(k)}$, namely,
	\begin{equation}\label{eq:2x2_rotation}
	\left(\!\begin{smallmatrix} \cos\phi_k & - e^{i\psi_k}\sin\phi_k \cr e^{-i\psi_k}\sin\phi_k & \cos\phi_k\end{smallmatrix}\right)
	\left(\!\begin{smallmatrix}{H}^{(k)}_{i_k i_k} & {H}^{(k)}_{i_k j_k}\cr\\[4pt] {H}^{(k)}_{j_k
		i_k} & {H}^{(k)}_{j_k j_k}\end{smallmatrix}\!\right) \left(\!\begin{smallmatrix}\cos\phi_k & e^{i\psi_k}\sin\phi_k \cr
	-e^{-i\psi_k}\sin\phi_k & \cos\phi_k\end{smallmatrix}\!\right) =
	\left(\!\begin{smallmatrix}{H}^{(k+1)}_{i_k i_k} & 0 \cr 0 & {H}^{(k+1)}_{j_k j_k}\end{smallmatrix}\!\right) .
	\end{equation}
	If the matrix $H$ is real, then $\psi_k\equiv 0$ and the transformation matrices are (real plane) Jacobi rotations.
	Unlike tridiagonalization-based methods, the Jacobi method does not preserve any zero structure.
	This method, originally  proposed by Jacobi for the real symmetric matrices \cite{jac-46} was rediscovered by Goldstine, Murray and von Neumann in \cite{gol-mur-neu-59}, and the extension to complex Hermitian matrices was done by Forsythe and Henrici \cite{for-hen-60}.  		 	
	An  instructive implementation with fine numerical details was provided by Rutishauser 	\cite{rutishauser-jacobi-66}, and an analysis of asymptotic convergence by Hari \cite{har-91-2}.
\end{itemize}
The beautiful simplicity of these methods  allows for quite some elegant generalizations. The Jacobi method, for instance, has been formulated and analyzed  in
 the context of Lie algebras \cite{Kleinsteuber-Helmke-Huper-Jacobi-CLA}, \cite{KLEINSTEUBER2009155}, and the \textsf{QR} method has its continuous form, the so-called Toda flow \cite{Isospectral-Watkins-1984}, \cite{QR-Toda-flow-Chu-1984}. 

\subsection{Backward stability in the conventional error analysis}\label{SS=backward-stability}
In finite precision (floating point) arithmetic, {not only are all processes  described above polluted by rounding errors, but
the iterations (\ref{eq:Hk-Lambda}) must be terminated at some appropriately chosen finite index $k_\star$}. In the $k$-th step, say, instead of $H^{(k)}$ we will have its computed approximation $\widetilde{H}^{(k)}$, for which a numerically unitary matrix $\widetilde{U}^{(k)}$  (i.e. $\| (\widetilde{U}^{(k)})^* \widetilde{U}^{(k)} - I\|_2 \leq O(n\roff)$) will be constructed, and the congruence transformation by $\widetilde{U}^{(k)}$ will be executed with rounding errors. By a backward error analysis, the new computed iterate $\widetilde{H}^{(k+1)}$ satisfies $\widetilde{H}^{(k+1)} = (\widetilde{U}^{(k)})^* ( \widetilde{H}^{(k)} + E_k)\widetilde{U}^{(k)}$. In general, this congruence is not a unitary similarity. Its software implementation, however, is carefully designed to ensure that both $E_k$ and  $\widetilde{H}^{(k+1)}$ are Hermitian.

\begin{figure}[ht]
	\setlength{\unitlength}{1.22cm}
	\centering
	\begin{picture}(10,3.0)(0,0)
	\put(0.7,2.5){$H$}
	\put(3.4,2.5){$\widetilde{H}^{(k_\star)} = \widehat{U}^* (H+\delta H)\widehat{U} = \widetilde{\Lambda} + \Omega(\widetilde{H}^{(k_\star)})$}
	\put(0.7,0.2){$H+\delta H$}
	\put(3.7,0.2){$\widehat{U}^* (H+\delta H - \widehat{U}\Omega(\widetilde{H}^{(k_\star)}) \widehat{U}^*)\widehat{U}=\widetilde{\Lambda} \approx \widetilde{U}^* H \widetilde{U}$}
	\put(1.1,2.55){\vector(0,-1){2.1}}
	\put(6.85,2.4){\vector(0,-1){1.7}}
	\put(7.0,1.7){set off-diagonal to zero}
	\put(7.0,1.3){equivalent to backward error} \put(7.0,0.9){$-\widehat{U}\Omega(\widetilde{H}^{(k_\star)})\widehat{U}^*$}
	\put(1.2,0.4){\vector(1,1){2.1}}
	\put(1.1,2.60){\vector(1,0){2.2}}
	\put(1.8,0.25){\vector(1,0){1.8}}
	\put(-0.25,1.8){backward}
	\put(-0.25,1.5){error $\delta H$}
	\put(2.2,1.2){exact similarity}
	\put(1.22,2.8){finite precision}
	\end{picture}
	\caption{\label{zd:FIG:eig_commutative_diagram} {Commutative diagram for a diagonalization routine in finite precision arithmetic. Since $\widetilde{U}$} is only numerically unitary, but close to an exactly unitary matrix $\widehat{U}$, we can claim that the computed output is close to an exact unitary diagonalization of a matrix close to the input $H$. This is sometimes called \emph{mixed stability}, i.e. both the input and the output must be changed to establish an exact relationship.}
\end{figure}
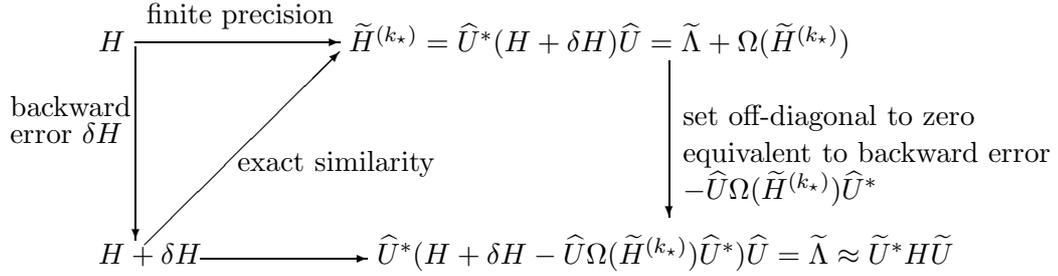

Since the computation has to be finite in time, one must also carefully determine the terminating index $k_\star$
such that the effectively computed matrix
\begin{equation}\label{eq:tlde-H-k}
\widetilde{H}^{(k_\star)} = (\widetilde{U}^{(k_\star-1)})^* ((\cdots ((\widetilde{U}^{(2)})^* ((\widetilde{U}^{(1)})^* (H + E_1) \widetilde{U}^{(1)} + E_2) \widetilde{U}^{(2)} + E_3) \cdots) +E_{k_\star-1})\widetilde{U}^{(k_\star-1)}
\end{equation}
is nearly diagonal, so that its sorted diagonal elements can be taken as approximate eigenvalues $\widetilde{\lambda}_1\geq\cdots\geq\widetilde{\lambda}_n$
of $H$. Let us assume that this sorting permutation is implicitly built in the transformation $\widetilde{U}^{(k_\star-1)}$, and let us write
$$
\widetilde{H}^{(k_\star)} = \widetilde{\Lambda} + \Omega(\widetilde{H}^{(k_\star)}),\;\;\mbox{where}\;\; \widetilde{\Lambda}=\begin{pmatrix} \widetilde{\lambda}_1 & 0& 0\cr
0 & \ddots & 0 \cr
0  &   0     & \widetilde{\lambda}_n\end{pmatrix},\;\; \widetilde{\lambda}_i = (\widetilde{H}^{(k_\star)})_{ii},\;\;i=1,\ldots, n,
$$
{and \,$\Omega(\cdot)$\, denotes the off-diagonal part of its matrix argument}.
If the eigenvectors are also needed, they are approximated by the columns of the accumulated product $\widetilde{U}$ of the transformations
$\widetilde{U}^{(k)}$. Since the accumulation is performed in finite precision, it can be represented as
$$
\widetilde{U} = (((( \widetilde{U}^{(1)} + F_1)\widetilde{U}^{(2)} + F_2) \widetilde{U}^{(3)} + \cdots ) + F_{k_\star-2})\widetilde{U}^{(k_\star-1)} .
$$
An error analysis ({see e.g. \cite[\S 6.5]{parlett-sevp-98}, \cite[Ch. 19]{book-higham}}) shows that for some small $\delta\widetilde{U}$ the matrix $\widehat{U}=\widetilde{U}+\delta\widetilde{U}$ is unitary, and
that there exists a Hermitian backward error $\delta H$ such that
\begin{equation}\label{zd:eq:tildeHk}
\widetilde{H}^{(k_\star)} = \widehat{U}^* (H+\delta H)\widehat{U}.
\end{equation}
Another tedious error analysis proves that $\|\delta H\|_2 \leq f(n)\roff \|H\|_2$, where the mildly growing function $f(n)$ depends on the details of each specific algorithm. The key ingredient in the analysis is the numerical unitarity of the transformations $\widetilde{U}^{(k)}$ in (\ref{eq:tlde-H-k}).
In the last step, the off-diagonal part $\Omega(\widetilde{H}^{(k_\star)})$ is deemed negligible and we use the approximate decomposition
$\widetilde{U}^* H \widetilde{U} \approx \widetilde{\Lambda}$. The whole procedure is represented by the diagram in Figure \ref{zd:FIG:eig_commutative_diagram},
and  summarized in Theorem \ref{zd:TM:eig_backward}.

%
\begin{theorem}\label{zd:TM:eig_backward}
	For each $i=1,\ldots, n$, let \,$\widetilde{\lambda}_{i}$ and $\widetilde{u}_i$ be the approximate eigenvalue and the corresponding approximate eigenvector of the Hermitian $n\times n$ matrix $H$, computed by one of the algorithms from  \S \ref{ss=Classical-Methods}. Then there exists a backward error $\Delta H$ and a unitary matrix \,$\widehat{U}$\, such that
	$H+\Delta H = \widehat{U}\widetilde\Lambda\widehat{U}^*$ and $\|\Delta H\|_2/\|H\|_2 \approx f(n)\roff$, $\|\widetilde U-\widehat U\|_2\approx O(n\roff)$, {where $\widetilde{U}$ is the $n\times n$ matrix with $i$-th column $\widetilde{u}_i$}. The finite precision computation can be represented by a commutative diagram as in Figure
	\ref{zd:FIG:eig_commutative_diagram}.
\end{theorem}
\noindent We have therefore a seemingly perfect situation: (i) The computed eigenvalues  $\widetilde{\lambda}_i$ are the exact eigenvalues of $H+\Delta H$, where $\|\Delta H\|_2/\|H\|_2$ is up to a factor of the dimension $n$ at the level of the roundoff unit $\roff$. (ii) By Theorem \ref{zd:TM:Weyl}, the absolute error in each $\widetilde{\lambda}_i$ is at most $\|\Delta H\|_2$.

\subsection{Case study: A numerical example}\label{zd:SSS:Example_3x3}
To put the framework of \S \ref{SS=backward-stability} under a stress test, we will compute the eigenvalues of a contrived $3\times 3$ real symmetric matrix, using the function \texttt{eig()} from the software package Matlab.  This function is based on the subroutine \texttt{DSYEV} from LAPACK \cite{LAPACK}, which implements a tridiagonalization-based \textsf{QR} algorithm.

\begin{example}\label{zd:EX:H3x3}
	{\em
	We use {Matlab R2010b} on a Linux workstation; the roundoff unit is $\roff \approx 2.2\cdot 10^{-16}$. The function ${\tt eig()}$ computes the approximate eigenvalues
	$\widetilde{\lambda}_1\geq\widetilde{\lambda}_2\geq\widetilde{\lambda}_3$
	($n=3$) of
	\begin{equation}\label{zd:eq:H3x3_eig()}
	H =
	\begin{pmatrix}
	10^{40}&  -2\cdot 10^{29} & 10^{19}\\[4pt]  -2\cdot 10^{29} & 10^{20} &10^9
	\\[4pt] 10^{19} & 10^{9} & 1
	\end{pmatrix}\;\;\mbox{as:}\;\;\;
	\begin{array}{r|c|}
	&   {\tt eig}(H) \\\hline
	\widetilde{\lambda}_1 & \;\;\; 1.000000000000000e+040\\
	\widetilde{\lambda}_2 &  -1.440001124147376e+020\\
	\widetilde{\lambda}_3 &  -1.265594217409065e+024
	\end{array}\;.
	\end{equation}
	Hence, the function ${\tt eig()}$ sees the matrix $H$ as indefinite with two negative eigenvalues.
	Following Theorem \ref{zd:TM:eig_backward} and Theorem \ref{zd:TM:Weyl},
	we know that $\max_{i=1:n}|\widetilde{\lambda}_i-\lambda_i| \leq O(\roff)\|H\|_2$, and that
	our computed eigenvalues are true eigenvalues of a nearby matrix $H+\delta H$, where $\|\delta H\|_2/\|H\|_2 \leq O(\roff)$.
	Since $n=3$, the effect of accumulated roundoff is negligible.

	However, to assess the quality of the approximation $\widetilde{\lambda}_i$ in terms of the number of its accurate digits, we need a bound to the relative error:
	\begin{equation}\label{zd:eq:H3x3_classical_bound}
	\frac{|\widetilde{\lambda}_i-\lambda_i|}{|\lambda_i|} \leq \frac{O(\roff)\|H\|_2}{|\lambda_i|}\leq \frac{O(\roff)\|H\|_2}{\min_{j=1:n}|\lambda_j|}
	\leq O(\roff)\|H\|_2 \|H^{-1}\|_2 \equiv O(\roff)\kappa_2(H).
	\end{equation}
	Thus, if $|\lambda_i| \approx \|H\|_2$, the computed approximation $\widetilde{\lambda}_i$ will have many correct digits.
	But if $|\lambda_i| \ll \|H\|_2$, then the above error bound cannot guarantee any correct digit in $\widetilde{\lambda}_i$.
	(It is immediately clear that $\lambda_1> 10^{40}$ and that $0<\lambda_3<1$, and thus $\kappa_2(H)>10^{40}$.)
	A second look at the matrix $H$ reveals its graded structure. In fact
	\begin{equation}\label{eq:H=DAD-3x3}
	H = D A D,\;\; A=\begin{pmatrix} 1 & -0.2\; & 0.1 \cr -0.2 & 1 & 0.1 \cr 0.1 & 0.1 & 1\end{pmatrix}, \; D = \begin{pmatrix} 10^{20} & 0 & 0 \cr
	0 & 10^{10} & 0 \cr 0 & 0 & 1\end{pmatrix},
	\end{equation}
	which means that $H$ is positive definite, and has a Cholesky factorization $H=LL^T$ (using Gershgorin circles one may immediately conclude that $A$ is positive definite.)
	
	Hence, the values of $\widetilde{\lambda}_2$ and $\widetilde{\lambda}_3$ in (\ref{zd:eq:H3x3_eig()}) are utterly wrong, although the result of
	${\tt eig}(H)$ is within the framework of \S \ref{SS=backward-stability}. Hence, the common routine of justifying the computed result by combining backward stability, ensured by unitary transformations, and the well-posedness in the sense of Weyl's Theorem is leading us to accept entirely wrong results. \hfill$\boxtimes$
}
\end{example}
	
\noindent At this point, one might be tempted to deem the matrix ill-conditioned and simply give up computing the tiniest eigenvalues (those with $|\lambda_i| < O(\roff) \|H\|_2$) to high accuracy, because they may not be well determined by the data; one may think it is simply not feasible.
{However, let
us explore this further:}

\begin{example}\label{zd:SSS:Example_3x3-b}
	{\em
	
	We continue with numerical experiments using the matrix $H$ from {Example \ref{zd:EX:H3x3}}; we apply ${\tt eig()}$ to the similar matrices
	$H_\zeta = P_\zeta^T H P_\zeta$, $H_\xi = P_\xi^T H P_\xi$, where $P_\zeta$, $P_\xi$ are the matrix representations of the permutations
	$\zeta=(3,2,1)$, $\xi=(2,1,3)$, respectively (in an application, this reordering could represent just another enumeration of
	the same set of variables and equations, thus describing precisely the same problem).
	Running ${\tt eig()}$ on these permuted matrices gives the following approximate spectra:
	\begin{equation}\label{zd:eq:H3x3_permuted}
	\begin{array}{r|c||c|}
	& {\tt eig}(H_\zeta)  &  {\tt eig}(H_\xi) \\ \hline
	\widetilde{\lambda}_1 &   1.000000000000000e+040 & 1.000000000000000e+040\\
	\widetilde{\lambda}_2 &   9.600000000000000e+019 & 9.600000000000000e+019\\
	\widetilde{\lambda}_3 &   9.750000000000001e-001 & 9.750000000000000e-001
	\end{array}\; .
	\end{equation}
	Similar values are computed with the permutation $\varpi=(3,1,2)$.
	With the permutations
	$(1,3,2)$ and $(2,3,1)$, however, ${\tt eig()}$ computes one single negative value, and the two positive
	values are identical to $\widetilde{\lambda}_1$ and $\widetilde{\lambda}_2$ in (\ref{zd:eq:H3x3_permuted}).
	All these computed approximate eigenvalues fit into the error estimate (\ref{zd:eq:H3x3_classical_bound}), but
	the qualitative difference is striking.
	
	Now let us compute, for the sake of experiment,  the eigenvalues of $H$ in two bizarre  ways: first as reciprocal values of the
	eigenvalues of $H^{-1}$, then as eigenvalues of numerically computed $H(H H^{-1})$.
	These are of course not very practical procedures, but we just want to see whether computing
	the eigenvalues of $H$ to high accuracy is warranted by the input. The results are as follows:
	\begin{equation}\label{zd:eq:H3x3_1/inv}
	\begin{array}{r|c||c|}
	& 1./{\tt eig}({\tt inv}(H))  &  {\tt eig}(H*(H\backslash H)) \\ \hline
	\widetilde{\lambda}_1 &   1.000000000000000e+040 & 1.000000000000000e+040\\
	\widetilde{\lambda}_2 &   9.600000000000000e+019 & 9.600000000000000e+019\\
	\widetilde{\lambda}_3 &   9.749999999999999e-001 & 9.750000000000002e-001
	\end{array}\; .
	\end{equation}
	Similar values to those in (\ref{zd:eq:H3x3_permuted}, \ref{zd:eq:H3x3_1/inv}) (up to a relative error of $O(\roff)$) are obtained either by calling
	${\tt eig}((H\backslash H)*H)$ or by using the result of ${\tt eig}(D^{-2},A)$ (here we use that the standard eigenproblem $Hx=\lambda x$ is equivalent to the generalized eigenproblem $A y=\lambda D^{-2}y$ with $y=Dx$ and $A$, $D$ as in (\ref{eq:H=DAD-3x3})).

	In view of all these numbers, what would be our best bet for the eigenvalues of $H$? Do the two smallest ones deserve to be computed better than already accepted and rationalized in Example \ref{zd:EX:H3x3}?
	
	Moreover, the function \texttt{chol()} computes the lower triangular Cholesky factor $L$ of $H$ as
	$$
	{\tt chol}(H)^T = \left(\begin{smallmatrix}
	\;\;\; 1.000000000000000e+020           &              0       &                  0 \\
	-2.000000000000000e+009  &  9.797958971132713e+009       &                  0 \\
	\;\;\;    9.999999999999999e-002  &  1.224744871391589e-001  &  9.874208829065749e-001
	\end{smallmatrix}\right) ,
	$$
	and the squared singular values of $L$  (the eigenvalues of $H$) are computed as
	$$
	{\tt svd}(L).^2
	=
	\left(\begin{matrix}
	1.000000000000000e+040, &
	9.600000000000002e+019, &
	9.750000000000000e-001 \end{matrix}\right) .
	$$
	Note that here the function ${\tt chol()}$, which is based on nonorthogonal transformations, correctly
	recognizes positive definiteness of $H$ and computes the triangular factor without difficulties. \hfill$\boxtimes$
}
\end{example}
Example \ref{zd:EX:H3x3} demonstrates that even in the symmetric case, computing the eigenvalues by the state-of-the-art software
tools may lead to difficulties and completely wrong results from a qualitative point of view. It is important to realize that,
in the framework described in \S \ref{SS=backward-stability}, the computed spectra (\ref{zd:eq:H3x3_eig()}), (\ref{zd:eq:H3x3_permuted})
 $H$ are all equally good and can be justified by a backward error analysis.
\begin{remark}
	{\em
	In Example \ref{zd:EX:H3x3} we specified that the results were obtained using Matlab R2010b on a Linux Workstation.
	On a Windows 10 based machine, the function ${\tt eig()}$ from Matlab R2015.a computes the eigenvalues of $H$ as
	$$
	{\tt eig}(H)= \left(\begin{matrix}
	1.000000000000000e+40, &
	9.900074641938021e-01, &
	-1.929211388222242e+23	
	 \end{matrix}\right) ,
	$$
	and e.g. $ 1./{\tt eig}({\tt inv}(H))$ returns the same result as in (\ref{zd:eq:H3x3_1/inv}). Numerical libraries (such as LAPACK \cite{LAPACK}, which is the computing engine for most of numerical linear algebra functions in Matlab) are often updated with improvements with respect to numerical robustness, optimizations with respect to run time etc. As a result, {{\em the same computation may return different results (better or worse) after a mere routine software update}}.\hfill$\boxtimes$
}
\end{remark}
Of course, if the initial $H$ is given with an uncertainty $\delta_0 H$ that is only known to be small in norm
($\|\delta_0 H\|_2/\|H\|_2 \ll 1$) then we cannot hope to determine the smallest eigenvalues in the case of
large condition number $\kappa_2(H)$. For instance, changing $H=\left(\begin{smallmatrix} 1 & 0 \cr 0 & \epsilon\end{smallmatrix}\right)$, $|\epsilon|\ll 1$, into $\left(\begin{smallmatrix} 1 & 0 \cr 0 & -\epsilon\end{smallmatrix}\right)$ is a small perturbation as measured in the operator norm $\|\cdot\|_2$, but it irreparably changes the smallest eigenvalue.

If, however, the data is given with smaller and more structured uncertainties, and if small matrix entries are not merely
noise, we ought to {do better/try harder}. Using the {customary}  norm-wise backward stability statement as a universal justification for errors in the result is not enough. This is best expressed by Kahan \cite{Kahan-Mindless}:
\emph{"The success of Backward Error-Analysis at explaining floating-point errors has been mistaken
	by many an Old Hand as an excuse to do and expect no better."}

\subsection{Computing the eigenvalues of positive definite matrices  with high relative accuracy}\label{SS=posdef-accurate}

Backward errors as analyzed in \S \ref{SS=backward-stability} are estimated in a matrix norm (usually $\|\cdot\|_2$ or $\|\cdot\|_F$). Unfortunately, a perturbation that is small in that sense may wipe out matrix entries that are in modulus much smaller than the matrix norm, and Example \ref{zd:EX:H3x3} shows that smallest eigenvalues may incur substantial damage. Demmel \cite{dem-92-2} showed that even for tridiagonal symmetric matrices, there are examples where tridiagonal \textsf{QR} with any reasonable shift strategy must fail to accurately compute the smallest eigenvalues.  However, some algorithms do produce better structured backward error that is gentler
to small entries, even if the initial matrix has no particular structure and its entries vary over several orders of magnitude. One consequence of more structured perturbation is that the condition number changes, and that a large standard condition number $\kappa_2(H)=\|H\|_2 \|H^{-1}\|_2$ does not necessarily imply that the computed eigenvalues will have large relative errors.

\subsubsection{Floating point perturbations and scaled condition numbers}\label{zd:SSS:RelativeFLoatingPointPert}
A closer look at some elementary factorizations, such as the Cholesky factorization of positive definite matrices, reveals that the standard norm-wise backward error analysis can be improved by estimating the relative errors in the individual entries.
We illustrate this kind of analysis with an important example of the Cholesky factorization of real symmetric
positive definite matrices.\footnote{Here the real case is cited for the sake of simplicity. An analogous analysis applies to the complex Hermitian case.}
\begin{theorem}\label{zd:TM:Demmel_on_Cholesky}(Demmel \cite{demmel-89-Cholesky})
	Let an $n\times n$ real symmetric matrix $H$ with strictly positive diagonal entries, stored in \textsf{IEEE} format with roundoff $\roff$, be input
	to the Cholesky algorithm. Let
	$$
	H=D H_s D,\; D={\rm diag}(\sqrt{H_{ii}})_{i=1}^n, \;\;
	( (H_s)_{ij}= \frac{H_{ij}}{\sqrt{H_{ii}H_{jj}}})
	$$
	and set  $\bfeta_C \equiv \frac{\max\{ 3, n\}\roff }{1- 2
		\max\{ 3, n\} \roff}>0$. Then:
	\begin{enumerate}
		\item If $\lambda_{\min}(H_s) > n \bfeta_C$, then the algorithm computes a
		lower triangular matrix $\widetilde{L}$ such that $\widetilde{L} \widetilde{L}^T=H+\delta H$, and for all
		$i,j=1,\ldots, n$ the backward error $\delta H$ can be bounded by
		$$
		|\delta H_{ij}| \leq \bfeta_C  \sqrt{H_{ii} H_{jj}}.
		$$
		Thus, {$H\!+\!\delta H \!\!=\!\! D(H_s\! +\! \delta
			H_s)D$,} {where $\delta H_s=D^{-1}\delta H D^{-1}$ satisfies}  $\max_{i,j}|(\delta H_s)_{ij}|\leq \bfeta_C\approx
			n\roff$.
		\item If $\lambda_{\min}(H_s) < \roff$, then there exists a sequence of simulated rounding errors that
		will cause the failure of the Cholesky algorithm.
		\item If $\lambda_{\min}(H_s) \leq - n \bfeta_C$, then the Cholesky algorithm will fail in floating point arithmetic with
		roundoff $\roff$.
	\end{enumerate}
\end{theorem}

Note that Theorem \ref{zd:TM:Demmel_on_Cholesky} does not assume that the matrix stored in the machine memory is positive definite.
Indeed, it can happen that we know a priori that our problem formulation delivers a positive definite
matrix, but the matrix actually stored in the computer memory is not definite, due to rounding errors.
The following simple example illustrates this.
\begin{example}\label{zd:EX:StifnessIndefinite}
	{\em
	(Cf. \cite[\S 11.1]{dgesvd-99})
	Consider the stiffness matrix of a mass spring system with $3$
	masses attached to a wall,
	with spring
	constants $k_1=k_3=1$, $k_2=\roff/2$:
	\begin{equation}\label{eq:mass-spring-3x3-H}
	{\boxminus\!\!\leftrightsquigarrow\!\!\blacksquare\!\!\leftrightsquigarrow\!\!
		\blacksquare\!\!\leftrightsquigarrow\!\!\blacksquare}\;\;\;\;
	H = \begin{pmatrix} k_1+k_2 & -k_2 & 0 \cr -k_2 & k_2+k_3 & -k_3\cr
	0 & -k_3 & k_3 \end{pmatrix}, \;\;{\lambda_{\min}(H)\approx
		\roff/4}.
	\end{equation}
	Here $\roff$ denotes the roundoff unit ({\tt eps}) in Matlab, so that $1+\roff/2$ is computed and stored as exactly $1$.
	The true and the computed assembled stiffness  matrix are, respectively,
	\begin{equation}\label{eq:mass-spring-3x3-tildeH}
	H = \begin{pmatrix} 1 + \frac{\roff}{2} &
	-\frac{\roff}{2} & 0 \cr -\frac{\roff}{2} & 1 +
	\frac{\roff}{2} & -1 \cr 0 & -1 & 1\end{pmatrix},\;\; \widetilde{H}
	=
	\begin{pmatrix}1 & -\frac{\roff}{2} & 0 \cr -\frac{\roff}{2} & 1 & -1 \cr 0 & -1 &
	1\end{pmatrix}.
	\end{equation}
	It is important to note here that the stored matrix $\widetilde{H}$ is component--wise close to $H$ with
	$${\displaystyle {|\widetilde{H}_{ij}-H_{ij}|}\leq \frac{\roff}{(2+\roff)} {|H_{ij}|}  < \frac{\roff}{2} {|H_{ij}|}}\;\;\mbox{for all}\;\;i,j.
	$$
	The matrix $H$ is by
	construction positive definite, whilst $\mathrm{det}(\widetilde{H})=-\roff^2/4$, and the smallest eigenvalue of $\widetilde{H}$ can be estimated as $\lambda_{\min}(\widetilde{H})\approx -\roff^2/8$.
	Hence, even computing the eigenvalues of $\widetilde{H}$ exactly could not provide any useful information
	about the smallest eigenvalue of $H$. The message of this example is: Once the data has been stored in the machine memory, the smallest eigenvalues can be so irreparably damaged that even exact computation cannot restore them.
	It is not hard to imagine how sensitive and fragile the computation of the smallest  eigenvalues of $H$ can be when
$\kappa_2(H)$ is large and the dimension of $H$ is  in the tens or hundreds of thousands, as e.g. in the case of discretizing
elliptic boundary value problems $\nabla\cdot (\bfa \nabla u)=f$ on $\Omega$, $u=g$ on the boundary of $\Omega$, where the scalar coefficient
field $\bfa$ on $\Omega$ varies over several orders of magnitude, see e.g. \cite{vavasis-96-FEMWC}. \hfill$\boxtimes$
}
\end{example}
\begin{remark}
	{\em
		Interestingly, if we consider assembling $H$ as in (\ref{eq:mass-spring-3x3-H}) as a mapping $(k_1,k_2,k_3)\mapsto H$, then the computation of $\widetilde{H}$ as in (\ref{eq:mass-spring-3x3-tildeH}) is not backward stable. There is no choice of stiffnesses $\widetilde{k}_1, \widetilde{k}_2$, $\widetilde{k}_3$ that would assemble (in exact arithmetic) to $\widetilde{H}$ that corresponds to three masses connected with springs as illustrated in (\ref{eq:mass-spring-3x3-H}). The reason is the indefiniteness of $\widetilde{H}$. On the other hand, the computation of $\widetilde{H}$ is perfectly forward stable. \hfill$\boxtimes$
	}
\end{remark}

\subsubsection{Characterization of well-behaved positive definite matrices}
\label{SSS=well-behaved-PD}
\noindent In some cases, the perturbation can be represented in the multiplicative form, i.e. $H+\delta H = (I+E) H (I+E)$ where $E$ can be bounded using the structure of $\delta H$. At the core of eigenvalue perturbation estimates is then Theorem \ref{zd:TM:Ostrowski} below.
For a detailed study, we refer to \cite{eis-ips-95} and \cite{Ren-Cang-Li-98-I}, \cite{Ren-Cang-Li-98-II}.
\begin{theorem} \label{zd:TM:Ostrowski}(Ostrowski \cite{ost-59})
	Let $H$ and $\widetilde{H}=Y^* H Y$ be Hermitian  matrices with eigenvalues
	$\lambda_1\geq\lambda_2\geq\cdots\geq\lambda_n$ and $\widetilde{\lambda}_1\geq\widetilde{\lambda}_2\geq \cdots\geq\widetilde{\lambda}_n$,
	respectively.
	Then, for all $i$,
	$$\widetilde{\lambda}_i =\lambda_i \xi_i,\;\;\quad \mbox{where }\quad
	\lambda_{\min}(Y^*\! Y)\leq \xi_i\leq \lambda_{\max}(Y^*\! Y).
	$$
\end{theorem}
The following theorem illustrates how the backward error of the structure as in Theorem \ref{zd:TM:Demmel_on_Cholesky}, combined with Theorem \ref{zd:TM:Ostrowski},  yields a sharp bound on the relative error in the eigenvalues.
\begin{theorem}\label{zd:TM:PertScaledPD}
	Let $\lambda_1\geq\cdots\geq\lambda_n$  and $\widetilde{\lambda}_1\geq\cdots\geq\widetilde{\lambda}_n$ be the eigenvalues of
$H=LL^T$ and of $\widetilde H=H+\delta H=\widetilde L\widetilde L^T$, respectively. If $\|L^{-1}\delta HL^{-*}\|_2 < 1$, then
	\begin{equation}\label{eq:dlambda}
	\max_i \left|\frac{\widetilde{\lambda}_i -\lambda_i}{\lambda_i}\right| \leq {\|
		H_s^{-1}\|_2} {\left\| \left[ \frac{\delta
			H_{ij}}{\sqrt{H_{ii}H_{jj}}}\right]_{i,j=1}^n \right\|_2},
	\end{equation}
{where $H_s$ is defined as in the statement of Theorem \ref{zd:TM:Demmel_on_Cholesky}}.
	(Recall the classical Weyl's theorem: $\max_{i}\left|\frac{\widetilde{\lambda}_i -\lambda_i}{\lambda_i}\right| \leq \kappa_2(H) \frac{\|\delta H\|_2}{\|H\|_2}$.)
\end{theorem}

\emph{Proof:}
Let $Y=\sqrt{I+L^{-1}\delta HL^{-*}}$. Then
$
H+\delta H = L(I+L^{-1}\delta HL^{-*})L^* = LYY^* L^*$ is similar to $Y^* L^* LY$, and we can equivalently compare	
the eigenvalues $\lambda_i(L^* L)=\lambda_i(H)$ and $\lambda_i(Y^* L^* LY)=\lambda_i(H+\delta H)$.
Now recall  Ostrowski's theorem: If $\widetilde M = Y^* M Y$, then, for all $i$, $\lambda_i(\widetilde M)=\lambda_i(M)\xi_i$, where
$\lambda_{\min}(Y^*\! Y)\leq \xi_i\leq \lambda_{\max}(Y^*\! Y)$.

Since $Y^*\! Y=I + L^{-1}\delta HL^{-*}$, we have
{$|\lambda_i(H)-\lambda_i(\widetilde H)|\leq \lambda_i(H) \|L^{-1}\delta H L^{-*}\|_2$}, with
\begin{eqnarray*}
	\|L^{-1}\delta H L^{-*}\|_2&=&\|L^{-1}D( {D^{-1}\delta H D^{-1}}) D L^{-*}\|_2=\|L^{-1} D ({\delta H_s}) D L^{-*}\|_2\\
	&\leq& \|L^{-1} D\|_2^2 \|\delta H_s\|_2=\|D L^{-*}L^{-1}D\|_2  \|{\delta H_s}\|_2\\
	&=& \| {(D^{-1} H D^{-1})^{-1}}\|_2 \|\delta H_s\|_2 = \| {H_s^{-1}}\|_2 \|\delta H_s\|_2,
\end{eqnarray*}	
{where we have denoted $\delta H_s=D^{-1}\delta H D^{-1}$ as in Theorem \ref{zd:TM:Demmel_on_Cholesky}}. The claim (\ref{eq:dlambda}) follows since ${\left(\delta H_s\right)_{ij}}=\delta H_{ij}/\sqrt{H_{ii}H_{jj}}$.
\hfill $\boxplus$

Since $\| H_s^{-1}\|_2 \leq \kappa_2(H_s)$, we see that, in essence, we have replaced the spectral condition $\kappa_2(H)$
with $\kappa_2(H_s)$, which behaves much better -- it is never much larger and it is potentially much smaller. In fact, $\|H_s^{-1}\|_2\leq \frac{n}{\|H_s\|_2}\min_{D=\mathrm{diag}}\kappa_2(D H D)$.
This claim is based on the following theorem.
\begin{theorem}\label{TM:VanDerSuis}
	(Van der Sluis \cite{slu-69})
	Let $H$ be a positive definite Hermitian matrix, $\Delta=\mathrm{diag}(\sqrt{H_{ii}})$ and $H_s = \Delta^{-1}H\Delta^{-1}$.
	Then
	$
	\kappa_2(H_s)\leq n \min_{D=diag}\kappa_2(D H D),
	$
{where the minimum is taken over all possible diagonal scalings $D$}.
\end{theorem}

\begin{example}\label{EX:cond-scond}{\em
	If we consider the matrix in Example \ref{zd:EX:H3x3}, we see that if $A\equiv H_s$, then $\|H_s^{-1}\|<1.4$ and  $\kappa_2(H_s) < 1.7$. This means
	that an algorithm with backward perturbation $\delta H$ of the form described in  Theorem \ref{zd:TM:PertScaledPD} may
	compute all eigenvalues of $H$ to nearly full machine precision (standard \textsf{IEEE} double precision with machine roundoff $\roff \approx 10^{-16}$) despite the fact that $\kappa_2(H)>10^{40}$.} \hfill$\boxtimes$
\end{example}
\begin{remark}\label{REM:scond-unit-i}{\em
{
The condition number $\kappa_2(H)$ is unitarily invariant: $\kappa_2(W^* H W)=\kappa_2(H)$ for any unitary matrix $W$. On the other hand, $\kappa_2((W^* H W)_s)$ can increase with a big factor. For instance, it is well known that there is a unitary $W$ such that $W^* H W$ has constant diagonal and thus (because of the homogeneity of the condition umber) $\kappa_2((W^* HW)_s)=\kappa_2(W^* H W)=\kappa_2(H)$, which can be much bigger that $\kappa_2(H_s)$, as illustrated in Example \ref{EX:cond-scond}. }	\hfill$\boxtimes$
}
\end{remark}
\noindent The number $\|H_s^{-1}\|_2$ can be interpreted geometrically in terms of the inverse distance to
singularity, as measured with respect to entry-wise perturbations. 
\begin{corollary} (Demmel \cite{demmel-89-Cholesky})
	Let $H=D H_s D$, where $D={\rm diag}(\sqrt{H_{ii}})_{i=1}^n$, and let
	$\lambda_{\min}(H_s)$ be the minimal eigenvalue of $H_s$.
	If $\delta H$ is a symmetric perturbation such that $H+\delta H$ is
	not positive definite, then
	$${\displaystyle \max_{1\leq i,j\leq
			n}\frac{|\delta H_{ij}|}{\sqrt{H_{ii}H_{jj}}} \geq
		\frac{\lambda_{\min}(H_s)}{n}=\frac{1}{n\|H_s^{-1}\|_2}}.$$
	If $\delta H = -\lambda_{\min}(H_s) D^2$, then ${\displaystyle
		\max_{i,j}\frac{|\delta
			H_{ij}|}{\sqrt{H_{ii}H_{jj}}}=\lambda_{\min}(H_s)}$ and $H+\delta H$
	is singular.
\end{corollary}
This means that in the case where $\|H_s^{-1}\|_2 > 1/\roff$, small entry-wise perturbations
can cause $H$ to lose the definiteness. Hence, if we assume no additional structure (such as
sparsity  pattern or signs distribution) a positive definite Hermitian matrix in floating point can be
considered numerically positive definite only if $\|H_s^{-1}\|_2$ is moderate (below $1/\roff$). The following
two results further fortify this statement.
\begin{theorem}\label{TM-SPD} (Veseli\'{c} and Slapni\v{c}ar \cite{ves-sla-93})
	Let $H=D H_s D$, {where $D={\rm diag}(\sqrt{H_{ii}})_{i=1}^n$}, be positive definite and let $c>0$ be a constant such that for
	all $\epsilon \in (0,1/c)$ and for all symmetric perturbations
	$\delta H$ with $|\delta H_{ij}|\leq\epsilon |H_{ij}|$, $1\leq
	i,j\leq n$, the ordered eigenvalues $\lambda_i$ and
	$\widetilde{\lambda}_i$ of $H$ and $H+\delta H$ satisfy ${\displaystyle
		\max_{1\leq i\leq n} \frac{|\widetilde{\lambda}_i - \lambda_i|}{\lambda_i}
		\leq c \epsilon.}$ Then ${\displaystyle \|H_s^{-1}\|_2 < (1+c)/2}$.
\end{theorem}

\begin{corollary}\label{COR-SPD}(Demmel and Veseli\'{c} \cite{dem-ves-92})
	Let $H$ be $n\times n$ positive definite and $\delta H = \eta D^2$, with any $\eta\in (0,\lambda_{\min}(H_s))$ and $D={\rm diag}(\sqrt{H_{ii}})_{i=1}^n$.
	Then for some index $\ell$ it holds that
	$$\frac{\widetilde{\lambda}_{\ell}}{\lambda_{\ell}} \geq
		 \sqrt[n]{1 + \eta \| H_s^{-1} \|_2}\equiv \sqrt[n]{1 + \max_{i,j}\frac{|\delta H_{ij}|}{\sqrt{H_{ii}H_{jj}}} \| H_s^{-1} \|_2} \approx 1 + \frac{\| H_s^{-1} \|_2}{n} \max_{i,j}\frac{|\delta H_{ij}|}{\sqrt{H_{ii}H_{jj}}} .
$$
\end{corollary}
Essentially,  if we have no additional structure (e.g. sparsity or sign pattern) accurate computation of all eigenvalues of Hermitian positive definite matrices in floating point is feasible if and only if $\|H_s^{-1}\|_2$ is moderate (as compared to $1/\roff$). In that case, allowing that the entries of $H$ are known up to small relative errors, computing the Cholesky factorization $H+\delta H=\widetilde{L}\widetilde{L}^*$ and working with $H$ in factored form via the computed factor $\widetilde{L}\approx L$ transforms the problem into the one of computing the \textsf{SVD} of the computed triangular factor. Since $\kappa_2(\widetilde{L})\approx \sqrt{\kappa_2(H)}$, the major part of the error is in the Cholesky factorization, as described in Theorem \ref{zd:TM:Demmel_on_Cholesky} and Theorem \ref{zd:TM:PertScaledPD}.

\subsubsection{Symmetric Jacobi algorithm for positive definite matrices}\label{SSS=symm-jac-pd}
Demmel and Veseli\'{c} \cite{dem-ves-92} proved that the symmetric Jacobi algorithm, when applied to the positive definite $H=H^T \in\mathbb{R}^{n\times n}$, produces backward errors that allow for direct application of Theorem \ref{zd:TM:PertScaledPD} at each iteration. To clarify, consider the commutative diagram of the entire process in Figure \ref{zd:FIG-symmj}.
\begin{figure}[H]
	\setlength{\unitlength}{5ex}
	\centering
	\begin{picture}(15,3.5)(0,0)
	\put(0,0){$\widetilde{H}^{(1)}+\delta\widetilde{H}^{(1)}$}
	\put(2.2,0.5){\vector(1,1){1.9}}
	\put(2.3,1.5){$\widehat{U}^{(1)}$}
	\put(2.3,2.75){$\widetilde{U}^{(1)}$}
	\put(3.05,0){$\widetilde{H}^{(2)}+\delta \widetilde{H}^{(2)}$}
	\put(5.3,0.5){\vector(1,1){1.9}}
	\put(5.4,1.5){$\widehat{U}^{(2)}$}
	\put(5.4,2.75){$\widetilde{U}^{(2)}$}
	\put(6.9,0){$\cdots$}
	\put(7.8,0){$\widetilde{H}^{(k_\star-1)}+\delta \widetilde{H}^{(k_\star-1)}$}
	\put(10.4,0.5){\vector(1,1){1.9}}
	\put(12.0,1.5){$\widehat{U}^{(k_\star-1)}$}
	\put(10.7,2.81){$\widetilde{U}^{(k_\star-1)}$}
	\put(0.1,2.5){$H=\widetilde{H}^{(1)}$}
	\put(1.1,2.35){\vector(0,-1){1.8}}
	\put(1.7,2.50){\vector(1,0){2.1}}
	\put(3.9,2.5){$\widetilde{H}^{(2)}$}
	\put(4.70,2.50){\vector(1,0){1.99}}
	\put(4.1,2.35){\vector(0,-1){1.8}}
	\put(6.9,2.5){$\cdots$}
	\put(7.6,2.50){\vector(1,0){1.3}}
	\put(9.1,2.5){$\widetilde{H}^{(k_\star-1)}$}
	\put(10.4,2.50){\vector(1,0){1.75}}
	\put(9.25,2.35){\vector(0,-1){1.8}}
	\put(12.2,2.5){$\widetilde{H}^{(k_\star)} = \widetilde{\Lambda} + \Omega(\widetilde{H}^{(k_\star)})$}
	\end{picture}
	\caption{\label{zd:FIG-symmj} The finite precision Jacobi algorithm for symmetric matrices. In a $k$-th step (the $k$-th commutative diagram),
		for the actually computed matrix 
		$\widetilde{H}^{(k+1)} \approx
		(\widetilde{U}^{(k)})^T \widetilde{H}^{(k)} \widetilde{U}^{(k)}$
		there exists a symmetric backward error $\delta \widetilde{H}^{(k)}$ and an orthogonal matrix $\widehat{U}^{(k)}$ such that		
		$\widetilde{H}^{(k+1)} = (\widehat{U}^{(k)})^{T}(\widetilde{H}^{(k)}+\delta\widetilde{H}^{(k)})\widehat{U}^{(k)}$ and
		$\| \widehat{U}^{(k)} -\widetilde{U}^{(k)}\|_2 = O(\roff)$. $\widetilde{\Lambda}$ is the diagonal part of $\widetilde{H}^{(k_{\star})}$. Compare with Figure \ref{zd:FIG:eig_commutative_diagram}.}
\end{figure}
\noindent Backward error analysis shows that for each iteration index $k$, 
\begin{equation}\label{eq:jac-rel-back-err}
|(\delta\widetilde{H}^{(k)})_{ij}|\leq O(\roff)\sqrt{(\widetilde{H}^{(k)})_{ii}(\widetilde{H}^{(k)})_{jj}}, \;\; 1\leq i,j\leq n, 
\end{equation}
which means that the relative error introduced in the $k$th step is governed by $\|(\widetilde{H}^{(k)})_s^{-1}\|_2$. (Here,  $(\widetilde{H}^{(k)})_s$ is defined analogously to $H_s$  in the statement of Theorem \ref{zd:TM:Demmel_on_Cholesky}, and we apply Theorem \ref{zd:TM:PertScaledPD}.)
The iterations are stopped at the first index $k_{\star}$ for which 
\begin{equation}
|(\widetilde{H}^{(k_{\star})})_{ij}|\leq O(\roff)\sqrt{(\widetilde{H}^{(k_{\star})})_{ii}(\widetilde{H}^{(k_{\star})})_{jj}}, \;\; 1\leq i,j\leq n,
\end{equation}
so that setting $(\widetilde{H}^{(k_{\star})})_{ij}$ to zero induces the perturbation $(\delta\widetilde{H}^{(k_{\star})})_{ij}=-(\widetilde{H}^{(k_{\star})})_{ij}$ of the type (\ref{eq:jac-rel-back-err}).
The overall accuracy depends on $\bfmu(H)=\max_{1\leq k\leq k_{\star}}\|(\widetilde{H}^{(k)})_s^{-1}\|_2$, which in practice is never much larger than $\|H_s^{-1}\|_2$. However, a formal proof of this remains an interesting open problem; for some discussion on this issue see \cite{mas-94}, \cite{drm-96-conbeh}.
An algorithm that computes the eigenvalues of any positive definite Hermitian $H$ to the accuracy determined by $\|H_s^{-1}\|_2$ (independent of $\bfmu(H)$) is given in \S \ref{SSS=Impl-Jacobi-eig}.

\subsection{Implicit representation of positive definite matrices}\label{SS=implicit-pd}

	In Example \ref{zd:EX:StifnessIndefinite}, $H\approx H_s$ and $\|H_s^{-1}\|_2\approx 1/O(\roff)$, and, by Theorem \ref{TM:VanDerSuis}, no diagonal scaling can substantially reduce its high condition number. One could argue that $H$ is ill-conditioned and that its smallest eigenvalue
	is not well determined by the data, i.e. by the matrix entries $H_{ij}$, and that it cannot be computed to any digit of accuracy. Indeed, the smallest eigenvalue has been lost at the very moment of storing the positive definite $H$ into the machine memory as the indefinite matrix $\widetilde{H}$, due to small relative changes (of the size of the machine roundoff unit) in the entries $H_{ij}$.  Hence, not even the exact computation with $\widetilde{H}$ could restore the information on the smallest eigenvalue of $H$.
	
	However, one may ask what data is actually given
	in this problem, and then argue that the data are the material properties (the stiffnesses $k_1, k_2, k_3$ of the springs) and the
	structure of the connections between the springs (adjacency).
	In fact, the stiffness matrix is usually assembled based on that information.
	In other words, $H$ can be written in a factored form
	as
	\begin{eqnarray}
		H &=& \begin{pmatrix} 1 & -1 & 0 \cr 0 & 1 & -1 \cr 0 & 0 & 1\end{pmatrix}
		\begin{pmatrix} k_1 & 0 & 0 \cr 0 & k_2 & 0 \cr 0 & 0 & k_3\end{pmatrix}
		\begin{pmatrix} 1 & 0 & 0 \cr -1 & 1 & 0 \cr 0 & -1 & 1 \end{pmatrix} \equiv B^T \mathrm{diag}(k_i)_{i=1}^3 B \label{eq:GTG1}\\
		&=& \begin{pmatrix}
			\sqrt{k_1} & -\sqrt{k_2} & 0 \cr 0 & \sqrt{k_2} & -\sqrt{k_3}\cr 0 & 0 & \sqrt{k3}\end{pmatrix}
		\begin{pmatrix} \sqrt{k_1} & 0 & 0 \cr -\sqrt{k_2} & \sqrt{k_2} & 0 \cr 0 & -\sqrt{k_3} & \sqrt{k_3}
		\end{pmatrix} = G^T G ,\;\; G = \mathrm{diag}(\sqrt{k_i})_{i=1}^3 B, \label{eq:GTG2}
	\end{eqnarray}
	thus clearly separating the adjacency from the material properties. {Furthermore, the \textsf{SVD} of the bidiagonal matrix $G$ can be computed to full machine precision \cite{dem-kah-90}, and once a high-accuracy \textsf{SVD} $G=U\Sigma V^T$ of $G$ is obtained, then $H = V \Sigma^2 V^T$ is the spectral decomposition of $H$.}
	
	If each $k_i$ is given with an initial uncertainty as $\widetilde{k}_i=k_i(1+\delta k_i/k_i)$, $\max_i|\delta k_i/k_i|\ll 1$, then in this factored
	representation we operate on
	$$
	\widetilde{G} = \begin{pmatrix} \sqrt{1+\delta k_1/k_1} & 0 & 0 \cr
	0 & \sqrt{1+\delta k_2/k_2} & 0 \cr 0 & 0 & \sqrt{1+\delta k_3/k_3}\end{pmatrix}
	\begin{pmatrix} \sqrt{k_1} & 0 & 0 \cr -\sqrt{k_2} & \sqrt{k_2} & 0 \cr 0 & -\sqrt{k_3} & \sqrt{k_3}
	\end{pmatrix} ,
	$$
	that is, $\widetilde{G}=(I+\Gamma)G$, where $\|\Gamma\|_2 \leq 0.5 \max_{i}|\delta k_i/k_i|$. By
	\cite{dem-kah-90} (see also the proof of Theorem \ref{zd:TM:PertScaledPD} and Theorem \ref{zd:TM:Eisenstat_Ipsen}) we know that the singular values
	of $G$ (and also the eigenvalues of $H$) are determined to nearly the same number of digits to which
	the coefficients $k_i$ are given, and that we can provably compute the singular values to that accuracy. This is in sharp contrast with the situation illustrated in Example \ref{zd:EX:StifnessIndefinite}.
	Hence, a different
	representation of the same problem is now perfectly well suited for numerical computations. The key for computing the eigenvalues of $H$ accurately is not to build $H$ at all, and to work directly with the parameters of the original problem.


	This example  raises an issue that is well known in the numerical linear
	algebra community, namely, that it is always advantageous to work with positive definite matrices {\em implicitly}. Since
	each positive definite matrix $H$ can be written as $H = A^* A$ with infinitely many choices for the full column rank
	(and in general rectangular) matrix $A$, we may find that in our specific situation such a factor is actually available. Here $A$ is not necessarily the Cholesky factor $L$, it need not be even square.
	Let us briefly comment a few well known examples.
	\begin{itemize}
		\item The solution of the linear least squares problem $\|Ax-b\|_2\longrightarrow \min$ with real full column rank $A$
		can be computed from the normal equations $A^T A x = A^T b$, but it is well known to numerical analysts that this is not a good idea because $H=A^T A$ satisfies $\kappa_2(H)=\kappa_2(A)^2$, see e.g. \cite[\S 2.1.4]{book-bjorck}. In other words if $A$ is $\epsilon=1/\kappa_2(A)$-close to singularity, then $H$ is $\epsilon^2$-close to some singular matrix. Notice that the positive definite matrix $H=A^T A$ is just an auxiliary object, not  part of the initial data ($A$, $b$), and that it has been invoked by the analytical characterization of the optimal $x$, where solving $Hx=b$ is considered simple since $H$ is positive definite.
		It turns out, in this case it is numerically advantageous to proceed by using the \textsf{QR} factorization with pivoting (or the \textsf{SVD}) of $A$, and
never to form and use $H$.
		\item Let $A$ be a Hurwitz-stable matrix {(i.e., all its eigenvalues lie on the left half-plane)} and suppose the matrix pair $(A,B)$ is controllable.  Computing the positive definite solution $H$ (controllability Gramian) of the Lyapunov matrix equation $A H + H A^* = -BB^*$ is difficult in the ill-conditioned cases and any
		numerical algorithm may fail to compute, in finite precision, the positive definite solution matrix $H$. Namely, if the solution $H$ is ill-conditioned, then it is close to the boundary of the cone of positive definite matrices and $H+\delta H$ may become semidefinite or indefinite, even for small forward error $\delta H$.	If the algorithm implicitly uses the assumed definiteness, it may fail to {run to completion/terminate} (analogously to the failure of the Cholesky decomposition of a matrix that is not numerically positive definite).
			
		It is better to solve the equation with the
		Cholesky factor $L$ of $H$ as the new unknown ($H=LL^*$, $L$ lower triangular with positive diagonal). Such approach was first advocated by Hammarling \cite{ham-82}, and later improved by Sorensen and Zhou \cite{Sorensen_Zhou:2003:DirectSylvLyap}. In this case, $L$ defines $H$ implicitly and $H=LL^*$ is positive definite. Also, since $\kappa_2(L)=\sqrt{\kappa_2(H)}$, the computed Cholesky factor $\widetilde{L}\approx L$ is more likely to be well-conditioned  and the implicitly defined solution $\widetilde{H}=\widetilde{L}\widetilde{L}^*$ is positive definite. Further, the Hankel singular values discussed in \S \ref{SS:scaling-Hankel-SVD}, can be computed directly from the Cholesky factors of the two Gramians, see \S \ref{S=PSVD}.
		\item In finite element computations, the symmetric positive definite $n\times n$ stiffness matrix $H$ can be assembled by factors -- the assembly process can be rewritten to produce a (generally, rectangular) matrix $A$ such that $H=A^T A$.
		This is the so called \emph{natural factor formulation}, see Argyris \cite{argyris-naturalFEM-75}. Such a formulation naturally leads to the Generalized Singular Value Decomposition (\textsf{GSVD}) introduced by Van Loan \cite{van-Loan-GSVD}. If e.g. $H_{ij}=\int_{a}^b \rho(x)\phi_i(x)\phi_j(x)dx$, $1\leq i,j\leq n$, and the integrals are evaluated by a quadrature formula $H_{ij}\approx \sum_{k=1}^m \omega_k \rho(x_k)\phi_i(x_k)\phi_j(x_k)$, then $H \approx A^T A$, where $A_{kj} = \sqrt{\omega_k \rho(h_k)}\phi_j(x_k)$, $1\leq k\leq m$, $1\leq j\leq n$. Note that $A = D \Phi$, where $D=\mathrm{diag}(\sqrt{\omega_k \rho(h_k)})_{k=1}^m$ and $\Phi_{kj}=\phi_j(x_k)$. Not only $\kappa_2(A)=\sqrt{\kappa_2(H)}$, but the most likely source of extreme ill-conditioning in $A$ (and thus in $H$) is clearly {isolated within} the diagonal matrix $D$. If an algorithm can exploit this and compute with a scaling invariant condition  number, then the essential and true condition number is that of $\Phi$, which depends on the choice of the basis functions $\phi_i(\cdot)$.
	\end{itemize}
	We conclude by noting that in all these examples, eigenvalue computations with $H$ can be equivalently done implicitly via the \textsf{SVD} or the \textsf{GSVD} of $A$.

\section{{Computing the} \textsf{SVD} {of arbitrary matrices}}\label{S=SVD}
Let $A$ be an $m\times n$ complex matrix. Without loss of generality, we assume $m\geq n$. The \textsf{SVD} $A=U\left(\begin{smallmatrix}\Sigma \cr 0\end{smallmatrix}\right)V^*$ implicitly provides the spectral decompositions of $A^* A$ and $AA^*$, and it is the tool of the trade in matrix computations and applications.  Numerical algorithms for computing the \textsf{SVD} are implicit formulations of  the diagonalization  methods for the Hermitian matrices. In the same way, the perturbation theory for the \textsf{SVD} is derived from the variational principles, and the error estimates of the computed singular values are derived from the combination of the backward error theory and the following classical result:

\begin{theorem}\label{TM:SVD-Weyl}
	Let the singular values of $A$ and $A+\delta A$ be
	$\sigma_1 \geq \cdots \geq\sigma_{\min(m,n)}$ and
	$\widetilde{\sigma}_1 \geq \cdots \geq\widetilde{\sigma}_{\min(m,n)}$, respectively. Then the distances between the corresponding singular values are estimated by a Weyl type bound
	$$
	\max_{i}|\widetilde{\sigma}_i-\sigma_i| \leq \| \delta A\|_2 .
	$$
	Further, the Wieland-Hoffman theorem yields
	${\displaystyle
	\sqrt{\sum_{i=1}^{\min(m,n)}|\widetilde{\sigma}_i-\sigma_i|^2}\leq \|\delta A\|_F .}
	$
\end{theorem}
To estimate the relative errors $|\widetilde{\sigma}_i - \sigma_i|/\sigma_i$, one uses perturbation in multiplicative form:
\begin{theorem}\label{zd:TM:Eisenstat_Ipsen} (Eisenstat and Ipsen, \cite{eis-ips-95})
	Let $\sigma_1\geq\cdots\geq\sigma_n$ and $\widetilde{\sigma}_1\geq\cdots\geq\widetilde{\sigma}_n$
	be the singular values of $A$ and $A+\delta A$, respectively. Assume that $A+\delta A$ can be
	written in the form of multiplicative perturbation $A+\delta A = \Xi_1 A \Xi_2$ and let
	$\xi=\max\{\|\Xi_1\Xi_1^T - I\|_2, \|\Xi_2^T\Xi_2-I\|_2\}$. Then
	$$
	|\widetilde{\sigma}_i - \sigma_i| \leq \xi \sigma_i,\;\;i=1,\ldots, n.
	$$
\end{theorem}
For  relative perturbation theory for the singular values and the singular vectors see \cite{Ren-Cang-Li-98-I}, \cite{Ren-Cang-Li-98-II}, and for an excellent review see \cite{Ipsen:1998:RPR}.

{We now describe and analyze three families of \textsf{SVD} algorithms, highlighting their different  properties regarding high-accuracy computations. For more details and further references see \cite[\S 8.6]{golub-vnl-4}, \cite[Ch. 58]{hogben14}.}

\subsection{Bidiagonalization-based methods}\label{SS=bidiagSVD}
Tridiagonalization (\ref{eq:tridiagonal}) of  $H=A^*A$ can be achieved implicitly by reducing $A$ to bidiagonal form \cite{gol-kah-65}
\begin{equation}\label{eq:bidiag}
U_1^* A V_1 = \begin{pmatrix} B \cr  0 \end{pmatrix} ,\;\;U_1, V_1 \;\;
\mbox{unitary},\;\;B=\left(\begin{smallmatrix} \alpha_1 & \beta_1 &  & \cr
& \alpha_2 & \ddots & & \cr
&  & \ddots & \beta_{n-1} \cr
& & & \alpha_n\end{smallmatrix}\right).
\end{equation}
The bidiagonalization process can be illustrated as follows {(as in \S \ref{ss=Classical-Methods} above,  $\star$ denotes an entry which has already
    been modified by the algorithm and set to its final value)}:
\begin{eqnarray*}
	A^{(1)} &=& U_{(1)}^T A =
	\left(\begin{smallmatrix}\star & \times & \times & \times\cr
		0 & \times & \times & \times\cr
		0 & \times & \times & \times\cr
		0 & \times & \times & \times\end{smallmatrix}\right),\;\;
	A^{(2)} = A^{(1)} V_{(1)} =
	\left(\begin{smallmatrix} \star & \star & 0 & 0\cr
		0 & \times & \times & \times \cr
		0 & \times & \times & \times \cr
		0 & \times & \times & \times \end{smallmatrix}\right),\;\;
	A^{(3)} =U_{(2)}^T A^{(2)} = \left(\begin{smallmatrix} \star & \star
		& 0 & 0 \cr
		0 & \star & {\times} & {\times}\cr
		0 & 0 & \times & \times \cr
		0 & 0 & \times & \times\end{smallmatrix}\right), \cr
	A^{(4)}&=&A^{(3)} V_{(2)}  = \left(\begin{smallmatrix} \star & \star
		& 0 & 0 \cr
		0 & \star & \star & 0\cr
		0 & 0 & \times & \times \cr
		0 & 0 & \times & \times\end{smallmatrix}\right),\;\;
	A^{(5)}= U_{(3)}^* A^{(4)}  = \left(\begin{smallmatrix} \star & \star
		& 0 & 0 \cr
		0 & \star & \star & 0\cr
		0 & 0 & \star & \star \cr
		0 & 0 & 0 & \star\end{smallmatrix}\right) = B = (U_{(3)}^* U_{(2)}^* U_{(1)}^*) A (V_{(1)} V_{(2)}) .
\end{eqnarray*}
Here $U_{(k)}$ and $V_{(k)}$ denote suitably constructed Householder reflectors, and their accumulated products form the matrices $U_1$ and $V_1$ in (\ref{eq:bidiag}). In the next step, the \textsf{SVD} of the bidiagonal $B = U_2 \Sigma V_2^*$,  can be computed by several efficient and elegant algorithms that implicitly work on the tridiagonal matrix $B^* B$; see e.g.  \cite{dem-kah-90}, \cite{mgu-eis-95}, \cite{GROSSER200345}. Combining the \textsf{SVD} of $B$ with the bidiagonalization (\ref{eq:bidiag}) yields the \textsf{SVD} of $A$:
$$
A = U_1 \begin{pmatrix} U_2 & 0 \cr 0 &  I \end{pmatrix}
\begin{pmatrix} \Sigma \cr  0 \end{pmatrix} (V_1 V_2)^* \equiv U \begin{pmatrix} \Sigma \cr  0 \end{pmatrix} V^* ,\;\;\Sigma = \left(\begin{smallmatrix} \sigma_1 &  \cr & \ddots & \cr & & \sigma_n\end{smallmatrix}\right),\;\;\sigma_1\geq\cdots\geq\sigma_n.
$$
Since only unitary transformations are involved,
we can prove existence of a backward error $\delta A$ and unitary matrices $\widehat{U}_1$, $\widehat{V}_1$ such that  the computed matrices $\widetilde{U}_1$, $\widetilde{V}_1$, $\widetilde{B}$ satisfy $\widetilde{U}_1\approx \widehat{U}_1$, $\widetilde{V}_1\approx\widehat{V}_1$ ($\widetilde{U}_1$, $\widetilde{V}_1$ are numerically unitary) and
\begin{equation}\label{eq:bidiagSVD}
A+\delta A = \widehat{U}_1\begin{pmatrix} \widetilde{B} \cr  0 \end{pmatrix}\widehat{V}_1^*  ,\;\;{\|\delta A\|_F} \leq \epsilon_1 {\|A\|_F}.
\end{equation}
It is important to know to what extent the singular values of $\widetilde{B}$ approximate the singular values of $A$. To that end, we invoke classical perturbation theory: if $\sigma_1(\widetilde{B})\geq\cdots\geq \sigma_{n}(\widetilde{B})$ are the singular values of $\widetilde{B}$ then applying Theorem \ref{TM:SVD-Weyl} to (\ref{eq:bidiagSVD}) yields
	\begin{equation}\label{eq:dsigma(B)}
	(i)\;\;\max_{i}|\sigma_i(\widetilde{B}) - \sigma_i| \leq \|\delta A\|_2;\;\;\;\;\;\;(ii)\;\;\sqrt{\sum_{i=1}^{n}|{\sigma}_i(\widetilde{B})-\sigma_i|^2}\leq \|\delta A\|_F \leq \epsilon_1 \|A\|_F.	
	\end{equation}
The computation of the \textsf{SVD} of $\widetilde{B}$ is also backward stable. If $\widetilde{\Sigma}$ is the diagonal matrix of the computed singular values, then, with some unitary matrices $\widehat{U}_2$, $\widehat{V}_2$ and some backward error $\delta \widetilde{B}$ we have
\begin{equation}\label{ex:SVD(B)}
\widetilde{B}+\delta\widetilde{B} = \widehat{U}_2\widetilde{\Sigma}\widehat{V}_2^* ,\;\; {\|\delta \widetilde{B}\|_F} \leq \epsilon_2 {\|\widetilde{B}\|_F}.
\end{equation}
 Both $\epsilon_1$ and $\epsilon_2$ are bounded by the roundoff $\roff$ times modestly growing functions of the dimensions.
The composite backward error of (\ref{eq:bidiagSVD}) and (\ref{ex:SVD(B)}) can therefore be written as
\begin{equation}\label{dA-composite}
A + \underbrace{\delta A + \widehat{U}_1 \begin{pmatrix} \delta\widetilde{B} \cr  0 \end{pmatrix}\widehat{V}_1^*}_{{\Delta A}} = \widehat{U}_1 \begin{pmatrix} \widehat{U}_2 & 0 \cr 0 &  I \end{pmatrix}
\begin{pmatrix} \widetilde{\Sigma} \cr  0 \end{pmatrix} (\widehat{V}_1\widehat{V}_2)^* \equiv \widehat{U}\begin{pmatrix} \widetilde{\Sigma} \cr  0 \end{pmatrix} \widehat{V}^* ,
\end{equation}
and the backward error is bounded in matrix norm as
$$
\|{\Delta A}\|_F \leq \epsilon_1 \|A\|_F + \epsilon_2 \|\widetilde{B}\|_F \leq
(\epsilon_1 + \epsilon_2 + \epsilon_1\epsilon_2) \|A\|_F.
$$
This is the general scheme of a bidiagonalization-based method. Depending on the method for computing the bidiagonal \textsf{SVD}, stronger statements are possible. For instance, if the \textsf{SVD} of $\widetilde{B}$ is computed with the zero-shift \textsf{QR} method \cite{dem-kah-90}, then all singular values of $\widetilde{B}$ (including the tiniest ones) can be computed to nearly full machine precision: if $\widetilde{\sigma}_1\geq\cdots\geq\widetilde{\sigma}_n$ are the computed values, then $|\widetilde{\sigma}_i - \sigma_i(\widetilde{B})| \leq O(n)\roff \sigma_i(\widetilde{B})$ for all $i$,  and the essential part of the error  $\widetilde{\sigma}_i - \sigma_i$ is committed in the bidiagonalization, so it is bounded in (\ref{eq:dsigma(B)}).
Note that, assuming  $A$ is of full rank and using (\ref{eq:bidiagSVD}),
\begin{equation}\label{dsigma-bidiag}
\max_i \frac{|\sigma_i(\widetilde{B})-\sigma_i|}{\sigma_i} \leq \frac{\|\delta A\|_2}{\sigma_{\min}} = \|A\|_2 \|A^\dagger\|_2 \frac{\|\delta A\|_2}{\|A\|_2} \equiv \kappa_2(A) \frac{\|\delta A\|_2}{\|A\|_2} \leq \sqrt{n}\epsilon_1 \kappa_2(A).
\end{equation}
Hence, although we can compute to nearly machine precision each, no matter how tiny, $\sigma_i(\widetilde{B})$, its value may be a poor approximation of the corresponding singular value $\sigma_i$ of $A$ if $\kappa_2(A)$ exceeds $O(1/\roff)$. For an illustration and explanation of how the reduction to bidiagonal form irreparably damages the smallest singular values see \cite[\S 5.3]{drm-xgesvdq}.  For an improvement of the backward error  (\ref{eq:bidiagSVD}) see \cite{bar-02-svd}.

\subsection{One-sided Jacobi \textsf{SVD}}\label{SS=One-sided-jacobi}
If the Jacobi method is applied to a real\footnote{Real matrices are used only for the sake of simplicity of the presentation.} symmetric positive definite matrix $H^{(1)}=H$, then the iterations $H^{(k+1)}=(V^{(k)})^T H^{(k)} V^{(k)}$ can be implemented implicitly: If one factorizes $H^{(k)}=(A^{(k)})^T A^{(k)}$, then $H^{(k+1)}=(A^{(k+1)})^T A^{(k+1)}$, where $A^{(k+1)} = A^{(k)} V^{(k)}$, $A^{(1)}=A$. If the pivot position at index $k$ is $(i_k,j_k)$, then the Jacobi rotation $V^{(k)}$ can be constructed from $A^{(k)}$ as follows: Let $d^{(k)}=(d_1^{(k)},\ldots,d_n^{(k)})$ be the diagonal of $(A^{(k)})^T A^{(k)}$. Compute
{${\displaystyle \xi_{i_k,j_k}=A^{(k)}(:,i_k)^T A^{(k)}(:,j_k)}$}, {where $A^{(k)}(:,s)$ denotes the $s$-th column  of $A^{(k)}$},  and
$${\vartheta_{i_k,j_k} =
\frac{d_{j_k}^{(k)}-d_{i_k}^{(k)}}{2\cdot \xi_{i_k,j_k}}}, \;\; t_k = {\frac{{\rm
sign}(\vartheta_{i_k,j_k})}{|\vartheta_{i_k,j_k}|+\sqrt{1+\vartheta_{i_k,j_k}^2}},\;\; c_k = \frac{1}{\sqrt{1+t_k^2}}}, \;\; s_k = t_k \cdot c_k.
$$
The transformation $A^{(k+1)}=A^{(k)} V^{(k)}$ leaves $A^{(k+1)}(:,\ell) = A^{(k)}(:,\ell)$ unchanged for $\ell\not\in\{i_k,j_k\}$, while
\begin{equation}\label{eq:one-sided-Jacobi}
\begin{pmatrix} A^{(k+1)}(:,i_k), & A^{(k+1)}(:,j_k)\end{pmatrix} =
			\begin{pmatrix} A^{(k)}(:,i_k), & A^{(k)}(:,j_k) \end{pmatrix} \begin{pmatrix}c_k & s_k \cr -s_k & c_k\end{pmatrix},
\end{equation}
and the squared column  norms are changed to
$
d_{i_k}^{(k+1)} = d_{i_k}^{(k)} - t_k \cdot \xi_{i_k,j_k}$,
$d_{j_k}^{(k+1)} = d_{j_k}^{(k)} + t_k\cdot\xi_{i_k,j_k}$.
If the accumulated product of the transformations $V^{(1)}\ldots V^{(k)}$ is needed, it can be updated analogously to (\ref{eq:one-sided-Jacobi}).
{Upon convergence, the limit of $(H^{(k)})_{k=1}^{\infty}$ is a diagonal positive definite matrix $\Lambda$, while the limit matrix of
$(A^{(k)})_{k=1}^{\infty}$ is $U\Sigma$, where the columns of
$U$ are orthonormal and $\Sigma = \sqrt{\Lambda}$. The columns of $U$ are the left singular vectors and the diagonal
matrix $\Sigma$ carries the singular values {of $A=U\Sigma V^T$, where
$V$, the accumulated product of Jacobi rotations, is orthogonal and has
the eigenvectors of $H$ as its columns}.}
This implicit application of the Jacobi method as an \textsf{SVD} algorithm is due to Hestenes \cite{hes-58}. An excellent implementation is provided by de Rijk \cite{rij-89}.

The key property of the Jacobi rotation, first identified by Demmel and Veseli\'{c} \cite{dem-ves-92}, is that the backward error in the finite precision implementation of (\ref{eq:one-sided-Jacobi}) is small in each pivot column, relative to that column. Hence, in a $k$th step, we have
\begin{equation}\label{eq:one-sided-Jacobi-backward}
\begin{pmatrix} \widetilde{A}^{(k+1)}(:,i_k), & \widetilde{A}^{(k+1)}(:,j_k)\end{pmatrix} =
\begin{pmatrix} \widetilde{A}^{(k)}(:,i_k) + \delta \widetilde{A}^{(k)}(:,i_k), & \widetilde{A}^{(k)}(:,j_k) + \delta \widetilde{A}^{(k)}(:,{j_k}) \end{pmatrix} \begin{pmatrix}\widetilde{c}_k & \widetilde{s}_k \cr -\widetilde{s}_k & \widetilde{c}_k\end{pmatrix},
\end{equation}
\begin{equation}\label{eq:one-sided-Jacobi-backward-columns}
\| \delta \widetilde{A}^{(k)}(:,i_k) \|_2 \leq \epsilon
\| \widetilde{A}^{(k)}(:,i_k)\|_2,\;\;
\| \delta \widetilde{A}^{(k)}(:,j_k) \|_2 \leq \epsilon
\| \widetilde{A}^{(k)}(:,j_k)\|_2 ,
\end{equation}
i.e.
$$
\widetilde{A}^{(k+1)} = (\widetilde{A}^{(k)} + \delta\widetilde{A}^{(k)}) \widetilde{V}^{(k)} = (I + \delta\widetilde{A}^{(k)} (\widetilde{A}^{(k)})^\dagger) \widetilde{A}^{(k)}\widehat{V}^{(k)}(I+E_k) ,
$$
where $\widehat{V}^{(k)}$ is orthogonal, $\|E_k\|_2\leq O(\roff)$.
By Theorem \ref{zd:TM:Eisenstat_Ipsen}, the essential part of the perturbation of the singular values, caused by $\delta \widetilde{A}^{(k)}$, is bounded by $\|\delta\widetilde{A}^{(k)} (\widetilde{A}^{(k)})^\dagger\|_2$ . Now let $D_k = \mathrm{diag}(\|\widetilde{A}^{(k)}(:,i)\|_2)$ and $\widetilde{A}^{(k)}_c = \widetilde{A}^{(k)} D_k^{-1}$. Then
$$
\|\delta\widetilde{A}^{(k)} (\widetilde{A}^{(k)})^\dagger\|_2 \leq
\| \delta\widetilde{A}^{(k)} D_k^{-1}\|_2 \| (\widetilde{A}^{(k)}_c)^\dagger\|_2 \leq \sqrt{2} \epsilon \| (\widetilde{A}^{(k)}_c)^\dagger\|_2 \leq \sqrt{2}\epsilon \kappa_2(\widetilde{A}^{(k)}_c).
$$
Note that $\widetilde{A}^{(k)}_c$ has unit columns and that, by Theorem \ref{TM:VanDerSuis}, $\kappa_2(\widetilde{A}^{(k)}_c)$ is up to a factor  $\sqrt{n}$ the minimal condition number over all diagonal scalings. The important property of the Jacobi algorithm, supported by overwhelming numerical evidence in \cite{dem-ves-92} is that $\max_{k\geq 1} \kappa_2(\widetilde{A}^{(k)}_c)$ is not much larger than $\kappa_2(A_c)$, where $A_c$ is obtained from $A$ by scaling its columns to unit Euclidean length. (Cf. \S \ref{SSS=symm-jac-pd}.)

Further, it is shown in \cite{drmac-97-rotations} that the Jacobi rotation can be implemented to compute the singular values in the full range of the floating point numbers. See \cite[\S 5.4]{drmac-HankelSVD-2015} for an example where the \textsf{SVD} is computed to high relative accuracy in \textsf{IEEE} double precision (64 bit) complex arithmetic despite the fact that $\sigma_{\max}/\sigma_{\min}\approx 10^{614}$.

Although more accurate than a bidiagonalization-based method, the one-sided Jacobi \textsf{SVD} has some drawbacks: its convergence may be slow, there is no sparsity structure to be preserved throughout the iterations  and each transformation is on a full dense matrix with low flop count per memory reference. These inconveniences can be alleviated by using the \textsf{QR} factorization as a preprocessor and a preconditioner for the one-sided Jacobi iterations.

\subsection{Jacobi \textsf{SVD} with \textsf{QR} preconditioning}\label{SS=J+QRCP}
In any \textsf{SVD} method, the \textsf{QR} factorization is a useful pre-processor, in particular in the case of tall and skinny matrices, i.e. $m\gg n$. Indeed, if $\Pi_r A\Pi_c= Q \begin{pmatrix} R \cr 0 \end{pmatrix}$ is the \textsf{QR} factorization with optional row and  column pivoting (encoded in the permutation matrices $\Pi_r$, $\Pi_c$), and $R=U_R \Sigma V_R^*$ is the \textsf{SVD} of $R$, then the \textsf{SVD} of $A$ is $A = \Pi_r^T Q \begin{pmatrix} U_R & 0 \cr 0 & I\end{pmatrix} \begin{pmatrix}\Sigma \cr 0\end{pmatrix} (\Pi_c V_R)^*$. In the sequel, we simplify the notation by assuming that the columns of the full column rank $A$ have been permuted so that $A\equiv \Pi_r A\Pi_c$.

If the one-sided Jacobi \textsf{SVD} is applied to $R$, then it implicitly diagonalizes $R^* R$. On the other hand, we can implicitly diagonalize $RR^*$ by applying the one-sided Jacobi to $R^*$. In that case the product of Jacobi rotations builds the matrix $U_R$. At first, there seems to be nothing substantial in this -- the \textsf{SVD} of $R$ and of $R^*$ are trivially connected. But, this seemingly simple modification is the key for faster convergence of the Jacobi iterations because $RR^*$ is more diagonally dominant than $R^* R$. There are deep reasons for this and the repeated \textsf{QR} factorization of the transposed triangular factor of the previous factorization is actually a simple way to approximate the \textsf{SVD},  \cite{mat-ste-93},
\cite{fer-par-95}, \cite{ste-97-qlp}.

The key is the column pivoting \cite{bus-gol-65} that ensures
\begin{equation}
|R_{ii}|\geq \sqrt{\sum_{k=i}^j
	|R_{kj}|^2},\;\;1\leq i\leq j \leq n.
\end{equation}
Such a pivoted \textsf{QR} factorization reveals the rank of $A$, it can be used to estimate the numerical rank, and it is at the core of many other methods, e.g. for the solution of least squares problems.

Let $A = A_c D_A$, $R = R_c D_c = D_r R_r$ with $D_A=\mathrm{diag}(\|A(:,i)\|_2)$, $D_c={\rm diag}(\| R(:,i)\|_2)$,
$D_r={\rm diag}(\| R(i,:)\|_2)$. Then $\kappa_2(A)=\kappa_2(R)$, $D_A = D_c$ and $\kappa_2(A_c)=\kappa_2(R_c)$, i.e. $R$ and $R_c$ inherit the condition numbers from $A$, $A_c$, respectively.
Furthermore, $R_r$ is expected to be better conditioned than $A_c$. It holds (see \cite{drmac-94-thesis,drm-99-Jacobi}) that $\kappa_2(R_r)$ is bounded by a function of $n$, independent of $A$, and that\footnote{Here the matrix absolute value is defined element-wise.}
\begin{equation}\label{eq:kappa-R-r}
\| R_r^{-1}\|_2 \leq \|\;|R_r^{-1}|\;\|_2 \leq\sqrt{n} \|\;|{R}_c^{-1}|\;\|_2 \leq n \| R_c^{-1}\|_2 .
\end{equation}
Hence, if $R$ can be written as $R = R_c D_c$ with well-conditioned $D_c$, then $R = D_r R_r$ with well conditioned $R_r$: $\| R_r^{-1}\|_2$ cannot be much bigger than $\| R_c^{-1}\|_2\equiv \|A_c^{\dagger}\|_2$, and it is potentially much smaller.

We now illustrate how an extremely simple (but carefully organized) backward error analysis yields sharp error bounds with a condition number that is potentially much smaller than the classical $\kappa_2(A)$.
The computed upper triangular $\widetilde{R}\approx R$ can be represented as the result of a backward perturbed \textsf{QR} factorization with an orthogonal matrix $\widehat{Q}$ and perturbation $\delta A$ such that
\begin{equation}\label{eq:qr:cols}
A + \delta A = \widehat{Q}\begin{pmatrix} \widetilde {R}\cr 0\end{pmatrix},\;\; \|\delta A(:,i)\|_2 \leq \roff_{qr} \|A(:,i)\|_2, \;\;i=1,\ldots, n.
\end{equation}
($\roff_{qr}$ is bounded by $\roff$ times a factor of the dimensions.)
If  we write this in the multiplicative form
\begin{equation}\label{eq:dA.A+}
A + \delta A = (I + \delta A A^\dagger) A,\;\; \|\delta A A^\dagger\|_2\leq \sqrt{n}\roff_{qr}\|A_c^\dagger\|_2 \leq
\sqrt{n}\roff_{qr}\kappa_2(A_c) ,
\end{equation}
and invoke Theorem \ref{zd:TM:Eisenstat_Ipsen}, we obtain
$$
\max_i \frac{|\sigma_i(A) - \sigma_i(\widetilde{R})|}{\sigma_i(A)}\leq 2 \|\delta A A^\dagger\|_2 + \|\delta A A^\dagger\|_2^2 \leq 2\sqrt{n} \roff_{qr} \kappa_2(A_c) + n (\roff_{qr} \kappa_2(A_c))^2.
$$
We conclude that the singular values of $\widetilde{R}$ are accurate approximations of the corresponding singular values of $A$, provided that $\kappa_2(A_c)$ is moderate. The key for invoking $\kappa_2(A_c)$ was (\ref{eq:qr:cols}), which was possible thanks to the fact that the \textsf{QR} factorization is computed by a sequence of orthogonal transformations that changed each column separately, without mixing them by linear combinations.\footnote{It is this mixing of large and small columns by orthogonal transformations oblivious of the difference in length that destroys the accuracy of the bidiagonalizaton. For illustrating examples see \cite{drm-xgesvdq}.}

If the one sided Jacobi \textsf{SVD} is  applied to $X=\widetilde{R}^T$, then its finite precision realization can be modeled as
\begin{equation}\label{eq:Jacobi-rows}
(X+\delta X)\widehat{V} \equiv (\widetilde{R} + \delta\widetilde{R})^T =  \widetilde{U}\widetilde{\Sigma},\;\;\|\delta X(i,:)\|_2 \leq \roff_J \|X(i,:)\|_2,\;\;i=1,\ldots, n,
\end{equation}
where $\widehat{V}$ s orthogonal and $\roff_J\leq O(n)\roff$. Note a subtlety here. The one sided Jacobi \textsf{SVD} is column oriented - the Jacobi rotations are designed to orthogonalize the columns of the initial matrix and in (\ref{eq:one-sided-Jacobi-backward}), (\ref{eq:one-sided-Jacobi-backward-columns}) the backward error analysis is performed column-wise. Here, for the purpose of the analysis, we consider the backward error row-wise.\footnote{We can also consider column-wise backward errors as in \S \ref{SS=One-sided-jacobi}, but this involves the behavior of scaled condition numbers of the iterates.} Hence, each row of $X=\widetilde{R}^T$, separately, has been transformed by a sequence of Jacobi rotations, and  we have (\ref{eq:Jacobi-rows}), where, in terms of the original variable $\widetilde{R}$, $\|\delta\widetilde{R}(:,i)\|_2 \leq \roff_J \|\widetilde{R}(:,i)\|_2 \leq \roff_J (1+\roff_{qr})\|A(:,i)\|_2$.

Finally, taking the \textsf{SVD} (\ref{eq:Jacobi-rows}) into (\ref{eq:qr:cols}), and writing numerically orthogonal $\widetilde{U}^T$ as $(I+E_u)^{-1}\widehat{U}^T$ with $\|E_u\|_2 \leq\roff_J$, we obtain
\begin{equation}
A + \underbrace{\delta A + \widehat{Q}\begin{pmatrix} \delta\widetilde{R}\cr 0\end{pmatrix}}_{\Delta A} = \widehat{Q} \begin{pmatrix} \widehat{V} & 0 \cr 0 & I_{m-n}\end{pmatrix} \begin{pmatrix} \widetilde{\Sigma}\cr 0 \end{pmatrix}(I+E_u)^{-1} \widehat{U}^T ,\;\;\|\Delta A(:,i)\|_2 \leq (\roff_{qr} + \roff_J(1+\roff_{qr}))\|A(:,i)\|_2 .
\end{equation}
This can be written as
\begin{equation}
\begin{pmatrix} \widehat{V}^T & 0 \cr 0 & I_{m-n}\end{pmatrix}\widehat{Q}^T (I+\Delta A A^\dagger) A\widehat{U}(I+E_u) = \begin{pmatrix} \widetilde{\Sigma}\cr 0 \end{pmatrix} ,
\end{equation}
where $\|\Delta A A^\dagger\|_2$ is estimated as in (\ref{eq:dA.A+}) with $\roff_{qr}+\roff_{J}(1+\roff_{qr})$ instead of $\roff_{qr}$,
and by Theorem \ref{zd:TM:Eisenstat_Ipsen},
\begin{equation}
\max_i \frac{|\sigma_i(A) - \widetilde{\Sigma}_{ii}|}{\sigma_i(A)}
\leq \max\{ 2\|\Delta A A^{\dagger}\|_2 + \|\Delta A A^{\dagger}\|_2^2, 2\|E_u\|_2 + \|E_u\|_2^2 \}
\end{equation}
\begin{theorem}\label{TM:jacobi-dsigma}
	Let $A$ be of full column rank and let its \textsf{SVD} be computed by Algorithm \ref{zd:ALG:eig:SVD:Jacobi} in finite precision with roundoff unit $\roff$, and let  $\widetilde{\sigma}_1\geq\cdots\geq\widetilde{\sigma}_n$ be the computed singular values. Assume no underflow nor overflow exceptions occur in the computation and let $\roff_{qr}$ and $\roff_J$ be as in (\ref{eq:qr:cols}) and (\ref{eq:Jacobi-rows}), respectively. Further, let $\roff_{\triangle}=\roff_{qr}+\roff_{J}(1+\roff_{qr})$.  Then
	\begin{equation}\label{dsigma-jacobi}
	\max_{i} \frac{|\widetilde{\sigma}_i - \sigma_i|}{\sigma_i} \leq 2\sqrt{n}\roff_{\triangle} \kappa_2(A_c) + n (\roff_{\triangle} \kappa_2(A_c))^2 .
	\end{equation}
\end{theorem}
\begin{algorithm}
	\caption{$(\Sigma, U, V) = \textsf{SVD}(A)$}
	\label{zd:ALG:eig:SVD:Jacobi}
	\begin{algorithmic}[1]
		%
		\STATE $(\Pi_r A)\Pi_c = Q \begin{pmatrix} R \cr 0 \end{pmatrix}$ \COMMENT{Rank revealing \textsf{QR} factorization; $\Pi_r$, $\Pi_c$  permutation matrices.}
		\STATE $X=R^T$; $X_\infty= X J_1 J_2 \cdots J_{\infty} = U_x \Sigma$ \COMMENT{One sided Jacobi \textsf{SVD}. }
		\STATE $V_x = J_1 J_2 \cdots J_{\infty}$
		\STATE $U = \Pi_r^T Q \begin{pmatrix} U_x & 0 \cr 0 & I\end{pmatrix}$, $V = \Pi_c Q_1 \begin{pmatrix} V_x & 0 \cr 0 & I\end{pmatrix}$
		\COMMENT{The \textsf{SVD} of $A$ is $A = U \begin{pmatrix}\Sigma \cr 0 \end{pmatrix} V^T$.}
	\end{algorithmic}
\end{algorithm}
Algorithm \ref{zd:ALG:eig:SVD:Jacobi} is the simplest form of the preconditioned one-sided Jacobi \textsf{SVD}. For a more sophisticated version, together with a detailed error analysis, including error bounds for the computed singular vectors, we refer to \cite{drm-ves-VW-1}, \cite{drm-ves-VW-2} and the LAPACK implementations \texttt{xGEJSV, xGESVJ}. The accuracy from Theorem \ref{TM:jacobi-dsigma} holds for any block oriented and parallelized implementation of Algorithm \ref{zd:ALG:eig:SVD:Jacobi}, see \cite{drm-block-jacobi-2010}.
Note that (\ref{dsigma-jacobi}) is preferred to the classical error bound (\ref{dsigma-bidiag}).
\begin{remark}
	{\em
The \textsf{QR} factorization  with column pivoting is at the core of many algorithms in a variety of software packages. Its first widely available robust implementation appeared in LINPACK \cite{LINPACK-UG} in 1979. It has been cloned and improved in LAPACK \cite{LAPACK}, and through LINPACK and LAPACK it has been incorporated into SLICOT, Matlab, and many other packages.
In 2008, it was discovered \cite{drmac-bujanovic-2008} that it contained a subtle instability that caused severe underestimation of the numerical rank of $A$ if $A$ is too close to the set of rank-deficient matrices. The problem was analyzed in detail and solved in \cite{drmac-bujanovic-2008}, and  the new code was incorporated into LAPACK in 2008, into SLICOT in 2010 (see \cite{buj-drm-slicot-2010}) and into ScaLAPACK in 2019 (see \cite{PXGEQPF}).
This is an example of how numerical instability can remain undetected for almost three decades, even in state-of-the-art software packages, inconspicuously producing bad results. This is also a warning and it calls for utmost rigor when developing and implementing numerical methods as scientific computing software. \hfill$\boxtimes$
}	
\end{remark}

\begin{remark}\label{REM:QR-xGESVD}
	{\em
Since the \textsf{QR} factorization is an efficient algorithm that reduces the iterative part to the $n\times n$ matrix $R$, the overall computation is more efficient in the case $m\gg n$. In fact, there is a crossover point for the ratio $m/n$ when even the bidiagonalization-based procedure is more efficient if it starts with the \textsf{QR} factorization and then bidiagonalizes $R$, see e.g. the driver subroutine \texttt{xGESVD} in LAPACK. Motivated by \cite{bar-02-svd}, we show in \cite{drm-xgesvdq} that, after using the \textsf{QR} factorization with pivoting as a preconditioner, the bidiagonalization becomes more accurate to the extent that in an extensive numerical testing the \textsf{QR} \textsf{SVD} (\texttt{xGESVD}) from LAPACK (applied to $R$ or $R^*$) matches the accuracy of the Jacobi method in Theorem \ref{TM:jacobi-dsigma}. This experimental observation seems difficult to prove. The algorithm is available in LAPACK as \texttt{xGESVDQ}.
	 \hfill$\boxtimes$
}
\end{remark}

\subsection{Accurate eigenvalues of positive definite matrices by the one-sided Jacobi algorithm}\label{SSS=Impl-Jacobi-eig}

From the discussion in \S \ref{zd:SSS:RelativeFLoatingPointPert}, it follows
that the Cholesky factorization is the perfect tool
for testing definiteness in floating point computation. Since our goal
is an accurate eigensolver for positive definite matrices with no
{\em a priori} given structure (e.g. zero or sign patterns or other
structural properties such as the Cauchy structure of the Hilbert matrix), we will use the Cholesky factorization to
test numerical positive definiteness.
%
%
It is remarkable that the following, very simple,
Algorithm \ref{zd:ALG:eig:CholJacobi},  proposed by
Veseli\'{c} and Hari \cite{ves-har-89}, provably achieves the optimal relative accuracy.
It is a combination of the Cholesky factorization and the one-sided Jacobi \textsf{SVD} algorithm.
\begin{algorithm}\label{ALG:jacobi-eig}
	\caption{$(\lambda,U) = \textsf{EIG}(H)$ ($H=H^T\in \mathbb{R}^{n\times n}$ positive definite)}
	\label{zd:ALG:eig:CholJacobi}
	\begin{algorithmic}
		%
		\STATE $P^T H P = L L^T$ \COMMENT{Cholesky factorization with pivoting.}
		\IF{$L$ computed successfully}
		\STATE $L_{\infty} = L \left< V \right>$
		\COMMENT{One--sided Jacobi \textsf{SVD} on $L$, without accumulation of the Jacobi rotations.}
		\STATE
		$\lambda_i = L_{\infty}(:,i)^T L_{\infty}(:,i), \;\;\;i=1,\ldots, n$ ;
		{$\lambda = (\lambda_1,\ldots,\lambda_n)$}.
		\STATE
		{$U(:,i) = {\displaystyle \frac{1}{\sqrt{\lambda_i}}} P L_{\infty}(:,i),\;\;\;i=1,\ldots,
			n$}.
		\ELSE
		\STATE Raise a warning flag: $H$ is not numerically positive definite.
		\STATE If the Cholesky factorization succeeded to compute $k$ columns of $L$, compute the \textsf{SVD} of the computed part $L(1:n,1:k)$
		(as above) and return $k$ positive eigenvalues with eigenvectors.
		\ENDIF
	\end{algorithmic}
\end{algorithm}

\noindent In the sequel, we will simplify the notation and assume that $H$
is already permuted, i.e. we replace $H$ with $P^T HP$ and analyze Algorithm \ref{zd:ALG:eig:CholJacobi}
with $P=I$. The following proposition is taken from \cite{drm-xgesvdq}.
\begin{proposition}\label{PR-BACKSEVP}
	Let $\widetilde L$, $\widetilde L_{\infty}$, $\widetilde U$,
	$\widetilde\lambda=(\widetilde\lambda_1,\ldots,\widetilde\lambda_n)$ be the
	computed approximations of $L$, $L_{\infty}$, $U$,
	$\lambda=(\lambda_1,\ldots,\lambda_n)$, respectively. Let
	$\widetilde\Lambda={\rm diag}(\widetilde\lambda_i)_{i=1}^n$. Then $\widetilde U
	\widetilde\Lambda \widetilde U^T=H+\Delta H$ with
	$$
	\max_{i,j}\frac{|\Delta H_{ij}|}{\sqrt{H_{ii}H_{jj}}} \leq
	\widetilde\bfeta_H\equiv \bfeta_C + (1+\bfeta_C)(2\roff_J+
	O(\roff)+O(\roff^2)).
	$$
\end{proposition}
{\em Proof}:\ We know that $\widetilde L \widetilde L^T = H+\delta
H\equiv\widetilde H$ with $|\delta H_{ij}|\leq \bfeta_C
\sqrt{H_{ii}H_{jj}}$ for all $i,j$. Further, we can write $\widetilde
L_{\infty}=(\widetilde L+\delta\widetilde L)\hat V$, where $\hat V$ is
orthogonal and $\|\delta\widetilde L(i,:)\|\leq\roff_J \|\widetilde
L(i,:)\|$ for all $i$. Let $\widetilde\Sigma={\rm
	diag}(\sqrt{\widetilde\lambda_1},\ldots,\sqrt{\widetilde\lambda_n})$.
A simple calculation shows that we can write $\widetilde U\widetilde\Sigma =
\widetilde L_{\infty}+\delta\widetilde L_{\infty}$, where $|\delta\widetilde
L_{\infty}|\leq \epsilon_\lambda |\widetilde L_{\infty}|$,
$0\leq\epsilon_\lambda\leq 3\roff$. Now it holds that $\widetilde U
\widetilde\Sigma^2\widetilde U^T = H + \delta H + E$, where for all $i,j$
$$
|E_{ij}| \leq 2 \left((\roff_J+\epsilon_\lambda
(1+\roff_J))+(\bfeta_J+\epsilon_\lambda (1+\roff_J))^2\right)
\sqrt{\widetilde H_{ii}\widetilde H_{jj}} \leq 2(\roff_J+O(\roff)+O(\roff^2))(1+\bfeta_C)\sqrt{H_{ii}H_{jj}}.
$$
\hfill $\boxtimes$

Strictly speaking, this proposition does not claim
backward stability of the eigendecomposition because $\widetilde U$ is
only nearly orthogonal; see \S  \ref{SS=backward-stability} and Figure \ref{zd:FIG:eig_commutative_diagram}. However, it is remarkable that $\widetilde
U\widetilde\Lambda\widetilde U^T$ recovers the original $H$ up to $O(n\roff)$
entry--wise relative errors in the sense of \S \ref{zd:SSS:RelativeFLoatingPointPert}.
As discussed in \cite{drm-xgesvdq}, in this situation one {suspects we} could use the \textsf{QR} algorithm-based \textsf{SVD} of $L_{\infty}$ to obtain equally good results, but a formal proof of high accuracy in this context is {still} lacking, see Remark \ref{REM:QR-xGESVD}.

\subsection{Eigenvalues of the pencil $HM-\lambda I$ and the \textsf{SVD} of matrix product}\label{S=PSVD}
In \S \ref{SS:scaling-Hankel-SVD}, we mentioned the importance of the eigenvalues of the product $HM$ (Hankel singular values), where $H$ and $M$ are real symmetric (or, more generally, Hermitian) positive definite matrices. If $H=L_h L_h^*$, $M=L_m L_m^*$ are the Cholesky factorizations of $H$ and $M$, then
$$
L_h^{-1}(HM)L_h = L_h^* L_m L_m^* L_h\equiv (L_m^* L_h)^* (L_m^* L_h) ,
$$
and the Hankel singular values are just the singular values of $A\equiv L_m^* L_h$. Set $D_h =\mathrm{diag}(\sqrt{H_{ii}})_{i=1}^n$, $H_s=D_h^{-1}HD_h^{-1}$, $L_{h,s}=D_h^{-1}L_h$, $D_m=\mathrm{diag}(\sqrt{M_{ii}})_{i=1}^n$, $M_s=D_m^{-1}MD_m^{-1}$, $L_{m,s}=D_m^{-1}L_m$. Note that $A=L_m^* L_h = L_{m,s}^* (D_m D_h) L_{h,s}$, where both $L_{h,s}$ and $L_{m,s}$ have rows of unit Euclidean length, and that $\kappa_2(L_{h,s})=\sqrt{\kappa_2(H_s)}$, $\kappa_2(L_{m,s})=\sqrt{\kappa_2(M_s)}$.
Based on our discussion in \S \ref{zd:SSS:RelativeFLoatingPointPert}, numerical positive definiteness of $H$, $M$ in the presence
of perturbations is feasible only if $\kappa_2(H_s)$ and $\kappa_2(M_s)$ are moderate; therefore we may assume that both
$L_{h,s}$ and $L_{m,s}$ are well conditioned.

This example motivates the study of numerical algorithms for computing the \textsf{SVD} of a matrix $A$ that is given in factored form $A = Z Y^*$, where  $Z\in\mathbb{C}^{m\times p}$ and $Y\in\mathbb{C}^{n\times p}$ are full column rank matrices such that
\begin{equation}\label{eq:zeta}
\bfzeta (Z,Y) \equiv \max\{ \min_{\Delta=\mathrm{diag}}\kappa_2(Z\Delta), \min_{\Delta=\mathrm{diag}}\kappa_2(Y\Delta)\}
\end{equation}
is moderate (below $1/\roff$).  Towards a more general situation, we may also write $Z Y^*$ as $X D Y^*$, where $D\in\mathbb{C}^{p\times p}$ is  diagonal, possibly very ill-conditioned.

\begin{example}
	{\em Computing the \textsf{SVD} of the product of matrices is an excellent example to illustrate the gap between the purely theoretical and actual computation in finite precision arithmetic: the simplest idea to compute the \textsf{SVD} of $Z Y^*$ (for given $Z$, $Y$) is to first compute the product $A=Z Y^*$ explicitly, and then reduce the problem to computing the \textsf{SVD} of $A$. The following example clearly illustrates the difficulty: If $\epsilon$ is such that $|\epsilon|<\roff$ (so $\pm 2+\epsilon$ is computed as $\pm 2$, and $\pm 1+\epsilon$ is computed as $\pm 1$ in finite precision) then
$$
\underbrace{\begin{pmatrix} 1 & \epsilon\cr -1 &
\epsilon\end{pmatrix}}_{Z} \begin{pmatrix} 2 & 2 \cr 2 & 1\end{pmatrix}=
\underbrace{\begin{pmatrix} 1 & 1\cr -1 &
	1\end{pmatrix}}_{X}
	\underbrace{\begin{pmatrix} 1 & 0 \cr 0 & \epsilon \end{pmatrix}}_{D}	
	\underbrace{\begin{pmatrix} 2 & 2 \cr 2 & 1\end{pmatrix}}_{Y^*} =
\begin{pmatrix} 2+2\epsilon & 2+\epsilon\cr -2+2\epsilon & -2+\epsilon\end{pmatrix}
$$
{will be computed and stored as}
$\widetilde{A}= \left(\begin{smallmatrix} 2 & 2 \cr -2 & -2\end{smallmatrix}\right)$, which means that the smallest singular value of order $|\epsilon|$ is irreparably lost. This problem is addressed by developing algorithms which avoid explicitly forming the matrix $A$. Instead, $Z$ and $Y$ are separately transformed in a sequence of iterations based on unitary matrices, see \cite{hea-lau-86}.
To ensure efficiency of the Kogbetliantz-type iteration, the matrices are unitarily transformed to triangular forms which are preserved throughout the iterations.
Since the entire computation relays on separate unitary transformations, the backward stability in the matrix norm is guaranteed. However, as illustrated in \S \ref{zd:SSS:Example_3x3}, this still does not guarantee high accuracy in the computed approximations of the smallest singular values.

To illustrate such a procedure and a numerical problem, the product $ZY^*$ is first transformed by
$(ZU_1^*)(U_1 Y^*)$ where $U_1$ is orthogonal such that $U_1 Y^*$ is upper triangular:
$$
U_1=\begin{pmatrix} \frac{1}{\sqrt{2}} &
\frac{1}{\sqrt{2}}\cr -\frac{1}{\sqrt{2}} &
\frac{1}{\sqrt{2}}\end{pmatrix},\;\;
U_1 Y^* =
\begin{pmatrix} \sqrt{8} & \frac{\sqrt{18}}{2}\cr 0
&-\frac{\sqrt{2}}{2}\end{pmatrix},\;\;
Z U_1^* =
\frac{1}{\sqrt{2}}\begin{pmatrix} 1+\epsilon & -1+\epsilon\cr
-1+\epsilon & 1+\epsilon\end{pmatrix}\approx
{\frac{1}{\sqrt{2}}\begin{pmatrix} 1
	& -1\cr -1 & 1\end{pmatrix}} .
$$
If $|\epsilon|$ is small relative to one,  $Z U_1^*$ will be computed and stored as an exactly singular matrix, and its smallest singular value will be lost in the very first step. It is  worth noticing that $Z$ has mutually orthogonal columns and that the computed version of $ZU_1^*$ is exactly singular, despite the fact that $U_1$ is orthogonal up to  machine precision. On the other hand, the value of $\bfzeta(Z,Y)$, defined in  (\ref{eq:zeta}), is easily seen to be less than $7$ in this example.
}
\end{example}

\subsubsection{An accurate algorithm}\label{SS=PSVD-Algorithm}

Applying the techniques from \S \ref{SS=J+QRCP}, we can easily construct an algorithm to compute the \textsf{SVD} of $A$ (given implicitly by $X$, $D$ and $Y$ as $A= X D Y^*$) with accuracy determined by $\bfzeta(X,Y)$ and independent of the condition number of $D$.
This allows for ill-conditioned $X$ and $Y$ as well, but such that ill-conditioning can be cured by diagonal scalings (i.e. moderate $\bfzeta(X,Y)$). Here we assume that $X$ and $Y$ are given either exactly, or that their columns are given up to small initial relative errors. Similarly, each diagonal entry of $D$ is given up to a small relative error.

\begin{algorithm}
	\caption{$(\Sigma, U, V) = \textsf{PSVD}(X, D, Y)$}
	\label{zd:ALG:PSVD}
	\begin{algorithmic}[1]
		%
		\STATE Factor $X = X_s \Delta_x$, where $\Delta_x=\mathrm{diag}(\|X(:,i)\|_2)_{i=1}^p$.
		Compute $Y_1 = Y D \Delta_x$.
		\STATE $Y_1 \Pi = Q \begin{pmatrix} R \cr 0 \end{pmatrix}$ \COMMENT{QR factorization with
			pivoting of $Y_1$.}
		\STATE $K=(X_s \Pi)R^*$ \COMMENT{Compute $K$ explicitly.}
		\STATE $K = U\begin{pmatrix}\Sigma \cr 0 \end{pmatrix} V_1^*$ \COMMENT{Compute the \textsf{SVD} of $K$ using the Jacobi method (Algorithm \ref{zd:ALG:eig:SVD:Jacobi}).}
		\STATE $V = Q \begin{pmatrix} V_1 & 0 \cr 0 & I\end{pmatrix}$
		\COMMENT{The \textsf{SVD} of $A$ is $A = U \begin{pmatrix}\Sigma \cr 0 \end{pmatrix} V^*$.}
	\end{algorithmic}
\end{algorithm}
\noindent To see why this algorithm is accurate (despite the fact that it uses the \textsf{SVD} of an explicitly computed matrix product, and that $D$ can be arbitrarily ill-conditioned) note the following:
\begin{itemize}
\item The column scaling in line 1. introduces entry-wise small relative errors, and it does not increase the condition number of the computation of the \textsf{QR} factorization in line 2. This is because the accuracy of the computed \textsf{QR} factorization of $Y_1$ is determined by $\min_{\Delta=\mathrm{diag}}\kappa_2(Y\Delta)$.
\item $R^*$ can be written as $R_r^* D_r^*$ with diagonal $D_r$ and well conditioned $R_r$. For Businger-Golub column pivoting \cite{bus-gol-65}, $\|R_r^{-1}\|_2$ can be bounded by $O(2^p)$ independent of $Y$, but if $Y$ is well conditioned, then 	$\|R_r^{-1}\|_2$ is expected to be
    at most $O(p)$. With so-called strong rank-revealing pivoting \cite{mgu-eis-96}, $\|R_r^{-1}\|_2$ can be bounded by $O(p^{1+(1/4)\log_2 p})$.
\item The matrix $K$ can be written as $K=K_c D_K$, where $D_K$ is diagonal and $K_c$ is well conditioned with equilibrated Euclidean column norms. The columns of $K$ are computed with small relative errors.
However, to preserve accuracy of even the tiniest singular values, the matrix multiplication must use the standard algorithm of cubic complexity. This is because the structure of the error of fast matrix
multiplication algorithms (e.g., Strassen) does not fit into the perturbation theory and cannot benefit from scaling invariant condition numbers.
\item In line 4., the Jacobi algorithm will compute the \textsf{SVD} with the accuracy determined by the condition number of $K_c$.
\end{itemize}
Hence, when it comes to computing the \textsf{SVD} with the condition number that is invariant under diagonal scalings, then we only need to carefully handle the scaling. The same argument applies to our claim that under the assumptions on the initial uncertainties in $X$, $D$ and $Y$, the \textsf{SVD} of $A\equiv XDY^*$ is determined to the accuracy with the condition number essentially given by
\begin{equation}\label{eq:mu}
\bfxi = \max\{ \|R_r^{-1}\|_2 \min_{\Delta=\mathrm{diag}}\kappa_2(X\Delta), \min_{\Delta=\mathrm{diag}}\kappa_2(Y\Delta)\}.
\end{equation}
For a more detailed analysis we refer the reader to \cite{drm-98-psvd}. For the case of more general $D$ see \cite{drm-98-tripletSVD}.

The decomposition of $A$ as $A = X D Y^*$, with diagonal $D$ and full column rank $X$ and such that $\min_{\Delta=\mathrm{diag}}\kappa_2(X\Delta)$ and  $\min_{\Delta=\mathrm{diag}}\kappa_2(Y\Delta)$ are moderate is called a \emph{rank-revealing} decomposition (\textsf{RRD}) of $A$. In the next section, we show that for some ill--conditioned matrices an accurate \textsf{RRD} can be computed to high accuracy that allows for computing accurate \textsf{SVD} by applying Algorithm  \ref{zd:ALG:PSVD}.

\subsection{Accurate \textsf{SVD} as \textsf{RRD}+\textsf{PSVD}}\label{S=RRD+PSVD}
Suppose we want to compute the \textsf{SVD} of $A$, but $A$ is so ill-conditioned that merely storing it in the machine memory may irreparably damage the \textsf{SVD}, or that all conventional  algorithms (cf. \S \ref{zd:SSS:Example_3x3}) fail due to extreme ill-conditioning (e.g. $A$ is the Hilbert or any other Cauchy or Vandermonde matrix).

An idea of how to try to circumvent such situation is presented in Example \ref{zd:EX:StifnessIndefinite}: the ill-conditioning of the matrix  is avoided by writing the matrix in factored form (\ref{eq:GTG1}, \ref{eq:GTG2}), using only a set of parameters $k_i$. The computed factored form is then used as input to an algorithm capable of exploiting the structure of the factors -- in this specific case, bidiagonal form. If we want to be able to tackle larger classes of difficult matrices, then we need to identify a factored form that is general enough and that we know how to use when computing the \textsf{SVD} to high accuracy, e.g. as with Algorithm \ref{zd:ALG:PSVD} in \S \ref{SS=PSVD-Algorithm}.
This is the basis of the approach introduced in \cite{dgesvd-99}. For more fundamental issues of finite precision (floating point) computation with guaranteed high accuracy see \cite{demmel_dumitriu_holtz_koev_2008}.

Suppose that $A$ can be written as $A = X D Y^*$, with $X$, $D$ and $Y$ as discussed in \S \ref{S=PSVD}, and that we have an algorithm that computes $\widetilde{X}=X+\delta X$, $\widetilde{D}=D+\delta D$, $\widetilde{Y}=Y+\delta Y$ such that\footnote{Alternatively, we may assume that $X$ and $Y$ are already well conditioned (thus properly scaled) and that the computed matrices satisfy
	$
	\|\delta X\|_2 \leq \epsilon_1 \|X\|_2,\;\;
	\|\delta Y\|_2 \leq \epsilon_2 \|Y\|_2,\;\;
	|\delta D_{ii}| \leq \epsilon_3 |D_{ii}|,\;\;i=1,\ldots, p.
	$}
\begin{equation}\label{eq:RRD-computed}
\|\delta X(:,i)\|_2 \leq \epsilon_1 \|X(:,i)\|_2,\;\;
\|\delta Y(:,i)\|_2 \leq \epsilon_2 \|Y(:,i)\|_2,\;\;
|\delta D_{ii}| \leq \epsilon_3 |D_{ii}|,\;\;i=1,\ldots, p.
\end{equation}
Write $\widetilde{D}$ as $\widetilde{D}=(I+E)D$, where $E$ is diagonal with $\|E\|_2\leq \epsilon_3$. Further, let $\Delta_X = \mathrm{diag}(\|X(:,i)\|_2)$, $X_c=X\Delta_X^{-1}$, $\delta X_c=\delta X\Delta_X^{-1}$; $\Delta_Y = \mathrm{diag}(\|Y(:,i)\|_2)$, $Y_c=Y\Delta_Y^{-1}$, $\delta Y_c=\delta Y\Delta_Y^{-1}$. Then
$$
\widetilde{X}\widetilde{D}\widetilde{Y}^* = (I+\delta_1 X X^\dagger) X D Y^*(I + \delta Y Y^\dagger)^*,\;\;\delta_1 X = \delta X + XE + \delta X E ,
$$
where the multiplicative error terms that determine the relative perturbations of the singular values  can be estimated as
\begin{eqnarray}
\|\delta_1 X X^{\dagger}\|_2 &\leq& \kappa_2(X_c) (\|\delta X_c\|_2 + \|E\|_2 + \|E\|_2\|\delta X_c\|_2) \leq \kappa_2(X_c) (\sqrt{p}\epsilon_1 + \epsilon_3 + \sqrt{p}\epsilon_1 \epsilon_3) \\
\|\delta Y Y^{\dagger}\|_2 &\leq&  \|\delta Y_c\|_2\|Y_c^\dagger\|_2 \leq \sqrt{p}\epsilon_2 \kappa_2(Y_c) .
\end{eqnarray}
Hence, if
$\bfzeta (X,Y)\equiv \max\{ \min_{\Delta=\mathrm{diag}}\kappa_2(X\Delta), \min_{\Delta=\mathrm{diag}}\kappa_2(Y\Delta)\}$
is moderate (below $1/\roff$), then the \textsf{SVD} of $A\equiv XDY^*$ can be accurately restored from the \textsf{SVD} decomposition of  $\widetilde{X}\widetilde{D}\widetilde{Y}^*$. For details see \cite{drm-98-psvd}, \cite{dgesvd-99}, \cite{Dopico-Moro-Mult-Error}.

The key advantages of  the factored representation are: \emph{(i)} The ill-conditioning is explicitly exposed in the ill-conditioned diagonal matrix $D$, and the factors $X$ and $Y$ are well-conditioned in the sense of (\ref{eq:zeta}). \emph{(ii)} The first errors committed in the computation are the small forward errors (\ref{eq:RRD-computed}) in $X$, $D$ and $Y$.

Hence, the problem is reduced to computing the decomposition $A=XDY^*$. This is solved on a case by case basis: first, a class of matrices is identified for which such a factorization is possible and then an algorithm for computing the decomposition $A= X DY^*$ is constructed. In the last step, the computed factors are given as input to Algorithm \ref{zd:ALG:PSVD}.

\subsubsection{\textsf{LDU}-based rank-revealing decompositions}\label{SS=svd-ldu-classes}
The \textsf{LDU} factorization with complete pivoting is used in \cite{dgesvd-99} as an excellent tool for providing \textsf{RRD}s of several important classes of matrices.
If $P_r A P_c = L D U$, with permutation matrices $P_r$, $P_c$, unit lower triangular $L$, diagonal $D$ and upper triangular $U$, then $X=P_r^T U$, $Y^* = U P_c^T$ yields $A = XDY^*$. Depending on the structure of $X$ and $Y$, we can deploy Algorithm \ref{zd:ALG:PSVD} (assuming only that $\bfzeta(X,Y)$ is moderate) or some other, more efficient, algorithm tailored for special classes of matrices. For instance, in  (\ref{eq:GTG1}, \ref{eq:GTG2}) the problem reduces to the \textsf{SVD} of a bidiagonal matrix and \textsf{QR} or \textsf{QD} algorithm can be applied. In some cases the sparsity pattern $\mathcal{S}$ (set of indices in the matrix that are allowed to be nonzero) and the sign distribution are the key properties for computing the singular values accurately.
We will here briefly mention few examples; for more detailed review see e.g. \cite[Ch. 59]{hogben14}, \cite{demmel_dumitriu_holtz_koev_2008}.

\paragraph{Acyclic matrices.}
 Let $A$ be  such that small relative changes of its nonzero entries (which are completely arbitrary, without any constraints) induce correspondingly small relative perturbations of its singular values (i.e. with the condition number $O(1)$). Then, equivalently, the associate bipartite graph $\mathcal{G}(A)$ is acyclic (forest of trees) and all singular values can be computed to high accuracy by a bisection method, see \cite{demmel-gragg-93}.  Bidiagonal matrices are acyclic and one can also use e.g. the zero-shift \textsf{QR} method \cite{dem-kah-90}. Also, the correspondence between the monomials in determinant expansion and perfect matchings in $\mathcal{G}(A)$ allows for accurate \textsf{LDU} factorization with pivoting.

\paragraph{Total sign compound (\textsf{TSC}) matrices.}
In some cases it is the sparsity and sign pattern $\mathcal{S}_{\pm}$ that facilitates an accurate \textsf{LDU} decomposition.
A sparsity and sign pattern $\mathcal{S}_{\pm}$ is \emph{total signed compound} (\textsf{TSC}) if
every square submatrix of every matrix $A$ with sign pattern $\mathcal{S}_{\pm}$ is either \emph{sign nonsingular} (nonsingular and
determinant expansion is the sum of monomials of
like sign) or \emph{sign singular} (determinant expansion
degenerates to sum of monomials, which are all zero).
Examples of \textsf{TSC} patterns are
	$$
	\left(\begin{smallmatrix}
	+ & + & 0 & 0 & 0 \\
	+ & - & + & 0 & 0 \\
	0 & + & + & + & 0 \\
	0 & 0 & + & - & + \\
	0 & 0 & 0 & + & +\end{smallmatrix}\right)
	,\;\;
	\left(\begin{smallmatrix}
	+ & + & + & + & + \\
	+ & - & 0 & 0 & 0 \\
	+ & 0 & - & 0 & 0 \\
	+ & 0 & 0 & - & 0 \\
	+ & 0 & 0 & 0 & -\end{smallmatrix}\right).
	$$
Suppose that every matrix $A$ with pattern
${\cal S}_{\pm}$ has the property that small relative changes of
its (nonzero) entries cause only small relative perturbations of
its singular values. Then  this property is equivalent with ${\cal S}_{\pm}$
being {total signed compound} (\textsf{TSC}). The \textsf{LDU} factorization with
complete pivoting $P_r A P_c = L D U $ of an \textsf{TSC} matrix $A$ can be computed so that all  entries of $L$, $D$, $U$ have small relative errors, and the framework of \S \ref{S=PSVD} applies. See \cite{dgesvd-99} for more details.

\paragraph{Diagonally scaled totally unimodular (\textsf{DSTU}) matrices.}
The $m\times n$ matrix $A$ is \emph{diagonally scaled
	totally unimodular (\textsf{DSTU})} if there exist diagonal matrices $D_1$, $D_2$ and a \emph{totally unimodular} $Z$ (all
minors of $Z$ are $-1$, $0$ or $1$) such that $A = D_1 Z D_2$.
To ensure that all entries of $L$, $D$ and $U$ are computed to high relative accuracy,  catastrophic cancellations (when subtracting intermediate results of the same sign) are avoided by predicting the exact zeros in the process of eliminations.  It can be shown that $\kappa_2(L)$ and  $\kappa_2(U)$
are at most $O(mn)$ and $O(n^2)$, respectively.

\paragraph{Cauchy matrices.}
\noindent Consider the \textsf{SVD} of a scaled (generalized) $m\times n$ Cauchy matrix
$$
C_{ij}=\frac{D_r(i) D_c(j)}{x_i + y_j},\;\; x, D_r \in\mathbb{R}^m,\;\;y, D_c\in\mathbb{R}^n.
$$
The key for the accuracy is in the fact that the \textsf{LDU} decomposition with full pivoting of $C$ can be computed
as a forward stable function of the vectors $x$ and $y$. More precisely, the decomposition $\Pi_1 C \Pi_2 = L D U$
($\Pi_1, \Pi_2$ permutations, $L$ unit lower triangular, $U$ unit upper triangular) is such that each entry is
computed to high relative accuracy and the triangular factors are well conditioned.
(In Algorithm \ref{zd:ALG:CauchyLDU}, the factorization is computed as $C=XDY^T\equiv (\Pi_1^T L) D (U\Pi_2^T)$.)
High accuracy of the computed factors follows from the fact that the Schur complement
can be recursively computed by explicit formulas involving only the initial vectors $x$ and $y$.
This is shown in Step 14. of Algorithm \ref{zd:ALG:CauchyLDU} by Demmel \cite{Demmel-99-AccurateSVD}.
The factors  can
be used in Algorithm \ref{zd:ALG:PSVD} as $X=L$, $Y^T = DU$, resulting in an accurate \textsf{SVD} of the product $LDU$.

\begin{algorithm}[hh]
	\caption{$(L,D,U,ir,ic) = $\textsf{CauchyLD}U$(x,y,D_r,D_c)$} \label{zd:ALG:CauchyLDU}
	\begin{algorithmic}[1]
		\STATE $m=max(size(x)); n=max(size(y)); p=min(m,n);$
		\FOR{$i = 1 : m$}
		\FOR{$j = 1 : n$}
		\STATE ${\displaystyle C(i,j) = \frac{ D_r(i) \cdot D_c(j)}{x(i)+y(j)}}$;
		\ENDFOR
		\ENDFOR
		\STATE $ir = [1:m]$; $ic = [1:n]$;
		\FOR{$k=1:p$}
		\STATE Find $(i_*,j_*)$ such that $|C(i_*,j_*)| = \max \{ |C(i,j)|\; :\; i=k,\ldots, m;\;j=k,\ldots, n \}$;
		\STATE ${\tt swap}(C(k,:),C(i_*,:))$;\ \   ${\tt swap}(C(:,k),C(:,j_*))$;\ \  ${\tt swap}(x(k),x(i_*))$;\ \   ${\tt swap}(y(k),y(j_*))$;
		\STATE ${\tt swap}(ir(k),ir(i_*))$;\ \  ${\tt swap}(ic(k),ic(j_*))$;
		\FOR{$i=k+1:m$}
		\FOR{$j=k+1:n$}
		\STATE ${\displaystyle C(i,j)=C(i,j) \frac{(x(i)-x(k))\cdot (y(j)-y(k))}{(x(k)+y(j))\cdot (x(i)+y(k))}}$;
		\ENDFOR
		\ENDFOR
		\ENDFOR
		\STATE $D = {\tt diag}(C)$ ;
		\STATE $X = {\tt tril}(G,-1) {\tt diag}(1./D) + {\tt eye}(m,n)$ ;
		\STATE $Y = ({\tt diag}(1./D)*{\tt triu}(G(1:n,1:n),1) + {\tt eye}(n))^T$ ;
		\STATE \COMMENT{$P=eye(m); \Pi_1=PP(:,ir)'; P=eye(n); \Pi_2 = P(:,ic); Y=\Pi_2*Y; X=\Pi_1^T*X;$}
	\end{algorithmic}
\end{algorithm}

\paragraph{Weakly diagonally dominant \textsf{M}-matrices.}
Suppose that the \textsf{M}-matrix $A=(A_{ij})\in\mathbb{R}^{n\times n}$ is diagonally dominant and that it is given with the off-diagonal entries $A_{ij}\leq 0$, $1\leq i\neq j\leq n$, and the row-sums $s_i = \sum_{j=1}^n A_{ij}\geq 0$. Note that this set of parameters determines the diagonal entries to high accuracy because $A_{ii}=S_i - \sum_{j, j\ne qi} A_{ij}$ has no subtractions/cancellations. Demmel and Koev \cite{demmel-koev-2004-M} showed that pivoted Gauss eliminations can be performed accurately in terms of the row sums and the off-diagonal entries, resulting in an accurate \textsf{LDU} decompositions, and an accurate \textsf{SVD}. For further details, see \cite{demmel-koev-2004-M}. With this unconventional matrix representation (off-diagonal entries and the row sums), it is possible to compute accurate \textsf{SVD} of diagonally dominant matrices, see \cite{svd-dd-matrix-Ye}, \cite{ldu-dd-matrix-dopico-koev}.
\begin{remark}{\em
	Due to pivoting, the factors $X$ and $Y$ are well conditioned. For example, if we factor the $100\times 100$
	Hilbert matrix $H_{100}$ (using a specialized version of Algorithm \ref{zd:ALG:CauchyLDU} for symmetric positive
	definite Cauchy matrices), as $H_{100} = X D X^T$ then $\kappa_2(X)\approx 72.24 \ll \kappa_2(H_{100}) > 10^{150}$.}
\end{remark}
\subsubsection{Con-eigenvalue problem for Cauchy matrices in the AAK theory}
More accurate numerical linear algebra impacts other approximation techniques in a variety of applications.
An excellent example is the case of $L^\infty$ rational approximations:
Haut and Beylkin \cite{haut-beylkin-coneig-2011} used Adamyan-Arov-Krein theory to
show that nearly $L^\infty$--optimal rational approximation on the unit circle   of
$
f(z) = \sum_{i=1}^n \frac{\alpha_i}{z-\gamma_i} + \sum_{i=1}^n \frac{\overline{\alpha_i}z}{1-\overline{\gamma_i}z}+\alpha_0
$
with a $m$-th order ($m<n$) rational function
$
r(z) = \sum_{i=1}^m \frac{\beta_i}{z-\eta_i} + \sum_{i=1}^m \frac{\overline{\beta_i}z}{1-\overline{\eta_i}z}+\alpha_0,\;\;\mbox{such that} \;\;\max_{|z|=1} , |f(z)-r(z)|\longrightarrow\min,
$
is numerically feasible if one can compute the con--eigenvalues and con--eigenvectors
of the positive definite generalized Cauchy matrix
$
{C=\left(\frac{\sqrt{\alpha_i}\sqrt{\overline{\alpha_j}}}{\gamma_i^{-1}-\overline{\gamma_j}}\right)} \in \mathbb{C}^{n\times n}.
$
In \cite{haut-beylkin-coneig-2011}
the con--eigenvalue problem $C u = \lambda \overline{u}$
is equivalently solved as the eigenvalue problem
$
\overline{C}C u = |\lambda|^2 u,
$
where $C$ is factored
as $C = X D^2 X^*$, and
$\overline{C}$ denotes the entry-wise complex conjugate matrix. The problem further reduces to
computing the \textsf{SVD} of the product $G = D X^T X D$, where $X$ is a complex matrix and $D$ is diagonal.
Such accurate rational approximation was successfully deployed in solving the initial boundary value
problem for the viscous Burger's equation {\cite{HAUT201383}}.

\subsubsection{Vandermonde matrices and the DFT trick}\label{SS=vandermonde}
In some cases, an \textsf{RRD} is not immediately available, but additional relations between structured matrices can be exploited.
 Demmel's algorithm for computing an accurate \textsf{SVD} of Vandermonde matrices \cite{Demmel-99-AccurateSVD} is a masterpiece
 of elegance. He used the fact that every $n\times n$ Vandermonde matrix $V=(x_i^{j-1})$ can be written as
$V = D_1 C D_2 F^*$, where $F$ is the unitary FFT matrix ($F_{ij}=\omega^{(i-1)(j-1)}/\sqrt{n}$, $\omega = \mathbf{e}^{2\pi\mathfrak{i}/n}$),
$D_1$ and $D_2$ are diagonal, and $C$ is a Cauchy matrix, i.e.,
\begin{equation}\label{zd:eq:VF-1}
(VF)_{ij} = \left[\frac{1-x_i^n}{\sqrt{n}}\right] \left[\frac{1}{\omega^{1-j}-x_i}\right] \left[\frac{1}{\omega^{j-1}}\right],\;1\leq i, j\leq n.
\end{equation}
After computing the \textsf{SVD} of the generalized Cauchy matrix $VF \equiv D_1 C D_2= U \Sigma W^*$, the \textsf{SVD} of $V$ is $V=U\Sigma (FW)^*$.
Note that in both cases the final step is the computation of the \textsf{SVD} of a product of matrices, based on
Algorithm \ref{zd:ALG:PSVD}.
This is turned into an accurate \textsf{SVD} of $\Vanderm_n(x)$, but with quite a few fine details, tuned to perfection in
\cite{Demmel-99-AccurateSVD}, \cite{demmel-koev-2006-V}. In particular, the possible singularity if some $x_i$ equals the floating point value of an $n$th root of unity is removable.
Demmel and Koev  \cite{demmel-koev-2006-V} extended this to polynomial Vandermonde matrices $V$  with entries $v_{ij}=P_i(x_j)$, where the $P_i$s are
orthonormal polynomials and the $x_j$s are the nodes.

\subsubsection{Toeplitz and Hankel matrices}
Let $\Hankel$ be a Hankel matrix, $\Hankel_{ij}=h_{i+j-1}$. The question is whether we can compute accurate singular values for any input vector $h$. This is equivalent to computing the singular values of the corresponding Toeplitz matrix $\mathcal{T}=P\Hankel$, where $P$ is the appropriate permutation matrix.
In general, a necessary condition to be able to compute all singular values of a square $A$ to high relative accuracy is that computing the determinant $\mathrm{det}(A)$ is possible to high accuracy. Applying this condition to the problem with Toeplitz matrices yields a negative result. It is impossible to devise an algorithm that can compute to high accuracy  the determinant of a Toeplitz or Hankel matrix for any input vector $h$.
In the fundamental work
\cite[\S 2.6]{Demmel-Dumitriu-Holtz}, it is shown that the obstacle in the complex case is the irreducibility of $\mathrm{det}(\mathcal{T})$ (over any field), and in the real case the problem is that $\nabla\mathrm{det}(\mathcal{T})$ has all nonzero entries on a Zariski open set.

However, in some settings the Hankel matrix $\Hankel$ is given implicitly as $\Hankel=\Vanderm^T D \Vanderm$, with suitable
Vandermonde $\Vanderm$ and diagonal $D$:
\begin{equation}\label{eq:H=VTDV}
\left(\begin{smallmatrix}
h_1 & h_2 & h_3 & \cdot & h_{n}\cr
h_2 & h_3 & \cdot & h_n & h_{n+1} \cr
h_3 & \cdot & \cdot & h_{n+1} & \cdot\cr
\cdot & h_n & h_{n+1} & \cdot & h_{2n-2}\cr
h_n & h_{n+1} & \cdot & h_{2n-2} & h_{2n-1}
\end{smallmatrix}\right) \!\! = \!\!
\left(\begin{smallmatrix}
1       & 1     & \cdot     & 1 & 1\cr
x_1     & x_2     & \cdot     & x_{n-1} & x_n \cr
x_1^2     & x_2^2   & \cdot     & x_{n-1}^2 & x_n^2\cr
\cdot    & \cdot & \cdot & \cdot & \cdot\cr
x_1^{n-1}     & x_2^{n-1}   & \cdot     & x_{n-1}^{n-1} & x_n^{n-1}
\end{smallmatrix}\right) \!\!\!
\left(\begin{smallmatrix}
d_1 &     &       &     &   \cr
& d_2 &       &     &    \cr
&     & \cdot &     &   \cr
&     &       & d_{n-1} &   \cr
&     &       &         &  d_n
\end{smallmatrix}\right) \!\!\!
\left(\begin{smallmatrix}
1     & x_1     & x_1^2     & \cdot & x_{1}^{n-1}\cr
1     & x_2     & x_2^2     & \cdot & x_{2}^{n-1} \cr
\cdot     & \cdot   & \cdot     & \cdot & \cdot\cr
1     & x_{n-1} & x_{n-1}^2 & \cdot & x_{n-1}^{n-1}\cr
1     & x_{n}   & x_n^2     & \cdot & x_n^{n-1}
\end{smallmatrix}\right) .
\end{equation}
If we refrain to compute $\Hankel$ (i.e. its vector $h$) explicitly and think of $\Hankel$ as parametrized by the numbers $x_i$, $d_i$, then accurate \textsf{SVD} of $\Hankel$ is possible. For tedious details and the full analysis we refer to \cite{drmac-HankelSVD-2015}.

\section{Computing accurate eigenvalues of Hermitian indefinite matrices}\label{S=HID}
The variational characterization of eigenvalues (Theorem \ref{TM-minmax}) and the resulting perturbation estimates (Theorem \ref{zd:TM:Weyl}) make no reference to the (in)definiteness (i.e. the inertia) of the Hermitian matrix $H$. Similarly, the state-of-the-art numerical software packages, such as LAPACK \cite{LAPACK}, use the generic routines for the Hermitian/symmetric eigenvalue problems that are backward stable in the sense of Theorem \ref{zd:TM:eig_backward} and accurate in the sense of (\ref{zd:eq:intro:rel_norm_error}). In \S \ref{S=HPD}, we discussed computation with high accuracy only for positive definite matrices (\S \ref{SS=posdef-accurate}).

When it comes to computing the eigenvalues with error bounds of the form (\ref{zd:eq:intro:rel_error}), there is a sharp distinction between  definite and indefinite matrices. 
For instance, for positive definite matrices the symmetric Jacobi eigenvalue algorithm is provably more accurate than the \textsf{QR} method \cite{dem-ves-92}, but in the case of indefinite matrices such general statement is not possible \cite{ste-95-qrbeatsjacobi}. Hence, new algorithms must be developed in order to guarantee reliable numerical results for indefinite matrices that are well-behaved with respect to finite precision diagonalization, in the sense that the computed eigenvalues satisfy the error bound (\ref{zd:eq:intro:rel_error}) with a moderate condition number $\bfkappa$.
Such matrices must be identified by the corresponding perturbation theory. In this section we give a brief review of theoretical results that have lead to good numerical algorithms.

\subsection{Perturbation theory for computations with indefinite matrices}\label{SS=HID-pert-theory}
Unfortunately, unlike the characterization of  positive definite matrices in \S \ref{SSS=well-behaved-PD}, perturbation theory provides no simple description of well-behaved  indefinite matrices.
The first important contribution to the theoretical understanding and algorithmic development was the analysis of the \emph{$\gamma$-scaled diagonally dominant} matrices \cite{bar-dem-90}, that are written as $H = D A D$, where $A=E+N$,  $E$ is diagonal with $E_{ii}=\pm 1$,  $D$ is diagonal with $D_{ii}=|H_{ii}|^{1/2}>0$, $N_{ii}=0$, and $\|N\|_2\leq \gamma <1$.
\begin{theorem} (Barlow and Demmel \cite{bar-dem-90}) Let $H = D A D$ be Hermitian $\gamma$-scaled diagonally dominant matrix with eigenvalues $\lambda_1\geq\cdots\geq\lambda_n$. Let $\delta H$ be a symmetric perturbation with $\| D^{-1}\delta H D^{-1}\|_2=\eta$, and let $H + \xi \delta H$ be $\gamma$-scaled diagonally dominant for all $\xi\in [0,1]$. If  $\widetilde{\lambda}_1\geq\cdots\geq\widetilde{\lambda}_n$ are the eigenvalues of $H+\delta H$, then, for all $i=1,\ldots, n$, 
$$
\frac{-\eta}{1-\gamma} + O(\eta^2) \approx
e^{-\eta/(1-\gamma)} - 1 \leq \frac{\widetilde{\lambda}_i-\lambda_i}{\lambda_i} \leq e^{\eta/(1-\gamma)} - 1 \approx \frac{\eta}{1-\gamma} + O(\eta^2) .
$$
\end{theorem}
Further, Barlow and Demmel \cite{bar-dem-90} showed that that a bisection algorithm can compute the eigenvalues of a $\gamma$-scaled diagonally dominant $H$  to high relative accuracy.  The key is that in this case, for any real $x$, the inertia of $H - x I$ can be computed with backward error $\delta H$ such that $\| D^{-1}\delta H D^{-1}\|_2$ is of the order of the machine precision. Recall, this means that the computed inertia is exact for the matrix $H+\delta H - x I$. Moreover, \cite{bar-dem-90} contains detailed analysis and computation of the eigenvectors, as well as extension of the results to symmetric $\gamma$-scaled diagonally dominant pencils $H-\lambda M$.

The seminal work of Barlow and Demmel initiated an intensive research, both for eigenvalue computations of Hermitian/symmetric matrices and the \textsf{SVD} of general and structured matrices. For the  Hermitian indefinite matrices, Veseli\'{c} and Slapni\v{c}ar \cite{ves-sla-93} generalized the results of \cite{bar-dem-90} to the matrices of the form $H = D A D$, $A = E+N$ with $E=E^*=E^{-1}$, $ED=DE$ and $\|N\|_2 <1$, and described an even larger class of well behaved matrices by identifying a new condition number,\footnote{The theory in \cite{ves-sla-93} has been developed for Hermitian pencils $H-\lambda M$ with positive definite $M$. Here we take $M=I$ for the sake of simplicity. } 
\begin{equation}
C(H) = \sup_{x\neq 0} \frac{|x|^T |H| |x|}{x^* \spabs H \spabs x},
\end{equation}
where $\spabs H \spabs = \sqrt{H^2}$ is the spectral absolute value of $H$, and $|H|$ is the element-wise absolute value, $|H|_{ij}=|H_{ij}|$. Note that $C(H)$ is finite for nonsingular $H$.

\begin{theorem} (Veseli\'{c} and Slapni\v{c}ar \cite{ves-sla-93}) \label{TM:VES-SLAP-T1}
Let $H$ and $H+\delta H$ be Hermitian with eigenvalues $\lambda_1\geq\cdots\geq\lambda_n$ and $\widetilde{\lambda}_1\geq\cdots\geq\widetilde{\lambda}_n$, respectively.
If the perturbation $\delta H$ is such that, for some $\eta<1$ and all $x\in\mathbb{C}^n$, {${\displaystyle |x^* \delta H x| \leq \eta x^* \spabs H \spabs x}$}, then $\widetilde{\lambda}_i=0$ if and only if $\lambda_i=0$, and for all nonzero $\lambda_i$'s
\begin{equation}\label{eq:eta-1}
\left| \frac{\widetilde{\lambda}_i - \lambda_i}{\lambda_i}\right| \leq \eta .
\end{equation}
\end{theorem}
The condition on $\delta H$ in this theorem is difficult  to check in practice, and in particular in case of the so-called floating point perturbations.\footnote{This term is used for typical errors occurring in the finite precision (computer) floating point arithmetic.}  The difficulty can be mitigated using $C(H)$ as follows. If we have $\delta H$ such that $|\delta H_{ij}|\leq \varepsilon |H_{ij}|$ for all $i, j$, then for any $x\in\mathbb{C}^n$, we have
$$
|x^* \delta H x| \leq |x|^T |\delta H| |x| \leq \varepsilon |x|^T |H| |x| \leq \varepsilon C(H) x^* \spabs H\spabs x,
$$
provided that $\spabs H \spabs$ is positive definite (i.e. $H$ nonsingular).

\begin{corollary}(Veseli\'{c} and Slapni\v{c}ar \cite{ves-sla-93})
	Assume that the matrix $H$ in Theorem \ref{TM:VES-SLAP-T1} is nonsingular, and that the Hermitian perturbation $\delta H$ is such that $|\delta H_{ij}|\leq \varepsilon |H_{ij}|$ for all $i, j$. If $\varepsilon C(H) <1$, then (\ref{eq:eta-1}) holds with $\eta = \varepsilon C(H)$.
\end{corollary}
Note that the condition $|\delta H_{ij}|\leq \varepsilon |H_{ij}|$ does not allow for perturbing zero entries; such strong condition cannot be satisfied in a numerical diagonalization process.
To get  more practical estimates, we must allow more general perturbations (see e.g. the conditions on $\delta H$ in Theorem \ref{zd:TM:Demmel_on_Cholesky} and Theorem \ref{zd:TM:PertScaledPD}), and have a more intuitive understanding of the factor $C(H)$.

To that end, $H$ is assumed nonsingular and $C(H)$ is estimated using the factored form $H = D A D$, where $D$ is diagonal matrix defined as the square root of the diagonal of $\spabs H \spabs$, $D = \mathrm{diag}(\spabs H\spabs)^{1/2}$. The role of the matrix $H_s$ from \S \ref{SSS=well-behaved-PD} has the matrix $\widehat{A} = D^{-1} \spabs H\spabs D^{-1}$.  
\begin{theorem}(Veseli\'{c} and Slapni\v{c}ar \cite{ves-sla-93})
Let 	$H = D A D$ be nonsingular, where $D$ is diagonal matrix defined as the square root of the diagonal of $\spabs H \spabs$, $D = \mathrm{diag}(\spabs H\spabs)^{1/2}$. Then 
$$
C(H)\leq \| |A| \|_2 \|\widehat{A}^{-1}\|_2 \leq \mathrm{Trace}(\widehat{A}) \| \widehat{A}^{-1}\|_2 \leq n \| \widehat{A}^{-1}\|_2 .
$$
If $\delta H$ is a Hermitian perturbation such that, for all $i, j$, $|\delta H_{ij}| \leq \varepsilon \sqrt{\spabs H\spabs_{ii} \spabs H\spabs_{jj}}$ and $\varepsilon n \| \widehat{A}^{-1}\|_2 < 1$, then (\ref{eq:eta-1}) holds for all eigenvalues, with $\eta =\varepsilon  n \| \widehat{A}^{-1}\|_2 $.
\end{theorem}	
Both definite and indefinite cases are nicely unified  by  Dopico, Moro and Molera \cite{dopico-moro-molera--2000-Weyl}, by showing that the argument used in the proof of Theorem \ref{zd:TM:PertScaledPD} extends, with careful application of the monotonicity principle, to indefinite matrices.
\begin{theorem}(Dopico, Moro and Molera \cite{dopico-moro-molera--2000-Weyl})
	Let $H$ and $H+\delta H$ be Hermitian with eigenvalues $\lambda_1\geq\cdots\geq\lambda_n$ and $\widetilde{\lambda}_1\geq\cdots\geq\widetilde{\lambda}_n$, respectively. Let $H$ be nonsingular, and let $H^{1/2}$ be any normal square root of $H$. If $\eta = \| H^{-1/2} \delta H H^{-1/2}\|_2 \leq 1$, then (\ref{eq:eta-1}) holds for all $i=1,\ldots, n$.
\end{theorem}

The difficulty in numerical computation is to have floating point backward errors that are compatible with the condition number. 
See e.g. the proof of Theorem \ref{zd:TM:PertScaledPD},  and note how in the relative error bound $\| L^{-1}\delta H L^{-*}\|_2$ the same diagonal scaling matrix results in small relative backward error and improved scaled condition number. At the same time, Algorithm \ref{ALG:jacobi-eig} has the backward errors of all iterations pushed back into the original matrix, and the structure of the error is compatible with this scheme, as shown in the proof of Proposition \ref{PR-BACKSEVP}. However, this is a special property of a particular algorithm. In general, it may be necessary to apply the perturbation estimate at each step in the sequence (\ref{eq:Hk-Lambda}), implemented as (\ref{eq:tlde-H-k}), and use the condition number of the current computed iterate $\widetilde{H}^{(k)}$.
It is important to note that the scaled condition numbers that are compatible with the structure of numerical errors are not invariant under unitary/orthogonal transformations, see Remark \ref{REM:scond-unit-i} and \cite{mas-94}, \cite{drm-96-conbeh}.
These issues have been successfully addressed in the algorithms that we review in \S \ref{SS=sym-indef-fact} and \S \ref{SS=OJ}. For simplicity, we consider only real symmetric matrices.
%

\subsection{Methods based on symmetric indefinite factorizations}\label{SS=sym-indef-fact}
The idea of using the pivoted Cholesky factorization of a positive definite $H$ to compute its spectral decomposition
via the \textsf{SVD} of its Cholesky factor (see \S \ref{SSS=Impl-Jacobi-eig}) can be also applied in the indefinite case. The first step is to obtain a symmetric indefinite factorization $H = G \mathcal{J} G^T$ with $\mathcal{J}$ diagonal, $\mathcal{J}_{ii} = \pm 1$,  $G$ well conditioned, and with backward error that allows for application of the perturbation theory with moderate condition numbers. This was first done by Slapni\v{c}ar \cite{slapnicar-GJGT-98}, who adapted the Bunch-Parlett factorization \cite{Bunch-Parlett}. An important feature of this factorization is that (as a result of pivoting), the matrix $G \mathrm{diag}(1/\|G(:,i)\|_2)$ is usually well conditioned, independent of the condition number of $H$.

For highly ill-conditioned matrices, accurate symmetric rank-revealing decompositions (\textsf{RRD}) $H= X D X^T$ have been computed in particular cases of structured matrices; see \cite{dopico-koev-2005-accurate} for symmetric totally nonnegative (\textsf{TN}) and Cauchy and Vandermonde matrices, and \cite{pelaes-moro-dstu-2006} for diagonally scaled totally unimodular (\textsf{DSTU}) and total signed compound (\textsf{TSC}) matrices.
Note that with $D = |D|^{1/2} \mathcal{J} |D|^{1/2}$, $\mathcal{J}_{ii}=\mathrm{sign}(D_{ii})$, and $G=X |D|^{1/2}$, the \textsf{RRD} $X D X^T$ can be written as $G\mathcal{J}G^T$.

In the next step, the eigenvalues and the eigenvectors of $H$ are computed using $G$ and $\mathcal{J}$ as the input matrices, i.e. $H$ is given implicitly by these factors ($H = G \mathcal{J} G^T$). 
We now briefly review two fundamentally different algorithms, which illustrate the development of accurate (in the sense of (\ref{zd:eq:intro:rel_error})) numerical methods for the symmetric indefinite eigenvalue problem.

\subsubsection{$J$-orthogonal Jacobi diagonalization}\label{SS=JJ}

Veseli\'{c} \cite{veselic-JS-Jacobi-93} noted that the factorization $H= G \mathcal{J} G^T$ can be used to compute the eigenvalues and eigenvectors of $H$ by diagonalizing the pencil $G^T G-\lambda \mathcal{J}$. Here $\mathcal{J}$ is a diagonal matrix with $\pm 1$ on its diagonal, and $G$ has $n$ columns and full column rank; in general $G$ can have more than $n$ rows.
The idea is to apply a variant of the one-sided Jacobi method, which we will now briefly describe.

In the $k$th step, $G^{(k+1)}=G^{(k)} V^{(k)}$ is computed from
$G^{(k)}$ using Jacobi plane rotations, exactly as in \S \ref{SS=One-sided-jacobi},
if $\mathcal{J}_{i_k i_k}$ and $\mathcal{J}_{j_k j_k}$ are of the same sign.
On the other hand, if $\mathcal{J}_{i_k i_k}$ and $\mathcal{J}_{j_k j_k}$ have opposite signs, then the Jacobi rotation is replaced with a hyperbolic
transformation
\begin{equation}\label{eq:hyp-rot}
\begin{pmatrix}
V^{(k)}_{i_k i_k} & V^{(k)}_{i_k j_k}\cr
V^{(k)}_{j_k i_k} & V^{(k)}_{j_k j_k}
\end{pmatrix} = 
\begin{pmatrix}
\cosh\zeta_k & \sinh\zeta_k \cr
\sinh\zeta_k & \cosh\zeta_k\end{pmatrix}, \quad \tanh 2\zeta_k = -\frac{2\xi_k}{d_{i_k}+d_{j_k}},
\end{equation}
$\xi_k=(G^{(k)})_{1:n,i_k}^T (G^{(k)})_{1:n,j_k}$,
$d_\ell=(G^{(k)})_{1:n,\ell}^T (G^{(k)})_{1:n,\ell}$, $\ell=i_k, j_k$.
The hyperbolic tangent is computed through
$${\displaystyle \tanh\zeta_k =
	\frac{{\rm sign}(\tanh 2\zeta_k)}{|\tanh 2\zeta_k|+\sqrt{\tanh^2
			2\zeta_k -1}}}.
$$
Note that $V^{(k)}$ belongs to the (unbounded) matrix group of  $\mathcal{J}$-orthogonal matrices, $(V^{(k)})^T \mathcal{J} V^{(k)}=\mathcal{J}$. (For basic properties of $\mathcal{J}$-orthogonal matrices, with applications in numerical linear algebra see \cite{higham-j-ort-review}.)
The limit of the $G^{(k)}$'s is $U\mathrm{diag}(\sigma_1,\ldots,\sigma_n)$; the $i$th column of $U$ is an eigenvector of $H$ associated with the eigenvalue $\lambda_i = \mathcal{J}_{ii}\sigma_i^2$.

Conceptually, this is an unusual approach, since, in the context of the symmetric eigenvalue problem, the use of orthogonal matrices
in a diagonalization process is considered natural and optimal, in particular in finite precision computation. Here, the elementary transformations matrices (\ref{eq:hyp-rot}) are orthogonal in the indefinite inner product induced by $\mathcal{J}$, and the matrix in the limit is $U\mathrm{diag}(\sigma_i)_{i=1}^n$ with $U^T U=I$.
The theoretical error bound for the computed eigenvalues contains a potential growth of the condition number, and in practice only a moderate growth has been observed, and the algorithm is considered accurate in the sense of (\ref{zd:eq:intro:rel_error}). For a detailed analysis and numerical evidence see \cite{slapnicar-thesis-92}, \cite{drm-har-93}, 
\cite{slapnicar-GJGT-98}, \cite{slapnicar-HAEVD-2002}, 
\cite{SLAPNICAR200057}. 
This approach can be applied to skew-symmetric problems $Sx=\lambda x$, $S\in\mathbb{R}^{2n\times 2n}$, see \cite{pie-93}.
In some applications, it is advantageous to formulate the problems in terms of the factors $G$ and $\mathcal{J}$, and  not to assemble the matrix $H$ at all, see \cite{veselic-GJGPT-2000}. 

\subsubsection{Implicit symmetric Jacobi method}\label{SSS=ISJM}
Dopico, Koev and Molera \cite{Dopico2009} 
 applied, analogously to the algorithm in \S \ref{SS=One-sided-jacobi}, the symmetric Jacobi method implicitly, i.e. by changing only the factor $G$ in $H=G\mathcal{J}G^T = X D X^T$. In the $X D X^T$ representation of the \textsf{RRD}, we can assume (by adjusting $D$) that $X$ has unit columns. Each matrix in the symmetric Jacobi algorithm is given implicitly as $H^{(k)} = G^{(k)} \mathcal{J} (G^{(k)})^T$, and one step of the method only computes $G^{(k+1)} = (V^{(k)})^T G^{(k)}$, thus implicitly defining $H^{(k+1)} = (V^{(k)})^T H^{(k)} V^{(k)} = G^{(k+1)} \mathcal{J} (G^{(k+1)})^T$. 
 
 This procedure can be preconditioned using the column pivoted \textsf{QR} factorization $ G \Pi = QR$, i.e. $Q^T H Q = R (\Pi^T \mathcal{J}\Pi) R^T$ and the implicit Jacobi scheme is applied to $R \widetilde{\mathcal{J}} R^T$, $\widetilde{\mathcal{J}} = \Pi^T \mathcal{J}\Pi$. As a result of preconditioning, the convergence may be substantially faster. 
 
This implicit Jacobi algorithm can be  implemented  to deliver the spectral decomposition with the accuracy  determined by $\kappa_2(X)$. This includes carefully designed stopping criterion, i.e. conditions to declare numerical convergence and to use the diagonal entries of the last implicitly computed $H^{(k)}$ as the approximate eigenvalues. For details, we refer to 
\cite{Dopico2009}.

\subsubsection{A remark on non-symmetric \textsf{TN} matrices}\label{SS=nonsymm-TN}
It has been noted in \cite{bar-dem-90}, \cite{ves-sla-93} that some of the perturbation estimates of the type (\ref{zd:eq:intro:rel_error}) extend to diagonalizable non-symmetric matrices. It is a challenging problem to identify classes of non-symmetric matrices for which eigenvalue computation with high accuracy is feasible, and to devise numerical algorithms capable of delivering such accuracy. Following the approach of \S \ref{SS=svd-ldu-classes}, the idea is to  use special matrix structure and the parameters that define its entries to find an implicit representation, in form of a decomposition, and then to apply a specially tailored algorithm.\footnote{Recall the example of the \textsf{SVD} of Vandermonde matrices in \S \ref{SS=vandermonde}.} 

The first successful breakthrough is the work of Koev \cite{koev-2005-TN} on totally nonnegative matrices. Totally nonnegative	(\textsf{TN}) matrices have all its minors nonnegative, and their eigenvalues are real and positive.   
A nonsingular \textsf{TN} matrix $A$ is characterized by a unique  bidiagonal decomposition $A=L_1\cdots L_{n-1} D U_{n-1}\cdots U_1$, where $D$ is diagonal, the $L_i$'s are unit lower bidiagonal, and the $U_i$'s are unit upper bidiagonal with additional zero and signs structure \cite{Gasca1996}. Such a decomposition follows from the Neville eliminations.	Koev \cite{koev-2005-TN} used this bidiagonal decomposition of a non-symmetric nonsingular \textsf{TN} matrix as the starting point for the first accurate algorithm, with detailed perturbation theory and error analysis, for computing eigenvalues of non-symmetric matrices. This led to a more general development of the numerical linear algebra of \textsf{TN} matrices, with new accurate algorithms for matrices derived from \textsf{TN} matrices \cite{koev-tn-2007}. 
Other examples of highly accurate solutions of non-symmetric eigenvalue problems include e.g. diagonally dominant \textsf{M}-matrices parametrized by the off-diagonal entries and the row sums \cite{eig-M-matrix-Ye}.

\subsection{Eigenvalue computation from the \textsf{SVD}}\label{SS=OJ}
 An accurate diagonalization  of an indefinite matrix can be derived from its accurate \textsf{SVD}
 decomposition, because the spectral and  the \textsf{SVD} decomposition are equal up to a multiplication with the inertia of $H$. Furthermore,  perturbation theory \cite{ves-sla-93} ensures that $H$ is well-behaved with respect to  computing eigenvalues if it is well-behaved with respect to computing the singular values. 
 To turn this into a robust numerical procedure, one must carefully recover the signs of the eigenvalues from the information carried by the singular vectors.  This has been done in \cite{dopico-molera-moro-2003-SEVP} with Algorithm \ref{zd:ALG:eig:indef:spanish} that turns any accurate \textsf{SVD} of a symmetric indefinite $H$ into an accurate spectral decomposition.

 \begin{algorithm}
 	\caption{$(\Lambda, Q) = \textsf{SVD2EIG}(H)$}
 	\label{zd:ALG:eig:indef:spanish}
 	\begin{algorithmic}[1]
 		%
 		\STATE $H = X D Y^T$. \COMMENT{Rank-revealing decomposition (\textsf{RRD}); not necessarily symmetric decomposition.}
 		\STATE $(\Sigma, Q, V) = \textsf{PSVD}(X,D,Y)$ \COMMENT{Algorithm \ref{zd:ALG:PSVD}.}
 		\STATE Recover the signs of the eigenvalues,
 		$\lambda_i=\pm\sigma_i$, using the structure of $V^T Q$.
 		\STATE Recover the eigenvector matrix $U$ using the structure of $V^T Q$.
 	\end{algorithmic}
 \end{algorithm}
\noindent Note that the first step aims at an accurate \textsf{RRD} and that symmetry of the decomposition is not the first priority. Also, for provable high accuracy for the computed eigenvalues and eigenvectors, the cases of multiple or tightly clustered singular values must be carefully analyzed; for further details, we refer to \cite{dopico-molera-moro-2003-SEVP}. Finally, note that we need to compute the full \textsf{SVD} of $H$, even if we only need  its eigenvalues.

\section{Acknowledgment}
The author wishes to thank Jesse Barlow, Jim Demmel, Froil\'{a}n  Mart\'{i}nez Dopico,  Vjeran Hari, Plamen Koev, Juan Manuel Molera Molera, Eberhard Pietzsch, Ivan Slapni\v{c}ar, Ninoslav Truhar, Kre\v{s}imir Veseli\'{c}, for numerous exciting discussions on accurate matrix computations, and in particular to Julio Moro Carre\~{n}o for encouragement to write this paper and for many useful comments that improved the presentation.

\bibliographystyle{plain}
\bibliography{eig-review_bib}
\end{document}